\documentclass[11pt]{article}
\usepackage{graphics}
\usepackage{eepic}
\usepackage{amsmath}
\usepackage{amscd}
\usepackage{amsfonts,latexsym,amssymb}
\input xy
\xyoption{all}
\xyoption{arc}
\xyoption{ps}
\xyoption{dvips}
\CompileMatrices

\newenvironment{proof}{\hspace{-\parindent}\textit{Proof.}}{\hfill $\Box$}

\newcommand{\SC}{{\mathcal{C}}}

\newcommand{\SE}{{\mathcal{E}}}
\newcommand{\SF}{{\mathcal{F}}}
\newcommand{\SG}{{\mathcal{G}}}
\newcommand{\SH}{{\mathcal{H}}}

\newcommand{\SM}{{\mathcal{M}}}

\newcommand{\SO}{{\mathcal{O}}}

\newcommand{\calS}{{\mathcal{S}}}

\newcommand{\SX}{{\mathcal{X}}}
\newcommand{\SY}{{\mathcal{Y}}}

\DeclareFontFamily{OT1}{rsfs}{}
\DeclareFontShape{OT1}{rsfs}{n}{it}{<->rsfs10}{}
\DeclareMathAlphabet{\curly}{OT1}{rsfs}{n}{it}

\newcommand{\CC}{\mathbb{C}}

\newcommand{\isom}{\cong}

\newcommand{\Hom}{\operatorname{Hom}}
\newcommand{\Quot}{\operatorname{Quot}}

\newcommand{\Spec}{\operatorname{Spec}}

\newcommand{\id}{\operatorname{id}}

\newcommand{\surj}{\twoheadrightarrow}

\newcommand{\too}{\longrightarrow}

\newcommand{\wt}{\widetilde}
\newcommand{\gitq}{{/\!\!/}}

\newcommand{\fc}{\mathfrak{C}}
\newcommand{\fcp}{\mathfrak{C'}}
\newcommand{\obc}{\operatorname{ob}\mathfrak{C}}

\newcommand{\Sch}{{Sch}}
\newcommand{\Sets}{{Sets}}
\newcommand{\Iso}{\operatorname{Iso}}
\newcommand{\groupoids}{{groupoids}}
\newcommand{\Bund}{\mathfrak{M}}
\newcommand{\FBund}{\underline{\Bund}}
\newcommand{\BBund}{{\mathcal{M}}}
\newcommand{\pgl}{{PGL(N)}}
\newcommand{\gl}{{GL(N)}}
\newcommand{\sqrtx}{\sqrt{x}}

\newtheorem{proposition}{Proposition}[section]
\newtheorem{theorem}[proposition]{Theorem}
\newtheorem{definition}[proposition]{Definition}
\newtheorem{lemma}[proposition]{Lemma}

\newtheorem{remark}[proposition]{Remark}
\newtheorem{example}[proposition]{Example}

\author{Tom\'as L. G\'omez \\
\\
Tata Institute of Fundamental Research\\
Homi Bhabha road, Mumbai 400 005 (India)\\
\texttt{tomas@math.tifr.res.in}}
\title{Algebraic stacks}
\date{9 November 1999}

\begin{document}

\maketitle

\begin{abstract}
This is an expository article on the theory of algebraic stacks.
After introducing the general theory, we concentrate in the example 
of the moduli
stack of vector budles, giving a detailed comparison with the
moduli scheme obtained via geometric invariant theory.
\end{abstract}

\section{Introduction}

The concept of algebraic stack is a generalization of the concept of
scheme, in the same sense that the concept of scheme is a generalization of
the concept of projective variety.
In many moduli problems, the functor that we want to study is not
representable by a scheme. In other words, there is no fine moduli
space. Usually this is because the objects that we want to parametrize
have automorphisms. But if we enlarge the category of schemes
(following ideas that go back to Grothendieck and Giraud, and were
developed by Deligne, Mumford and Artin) and consider algebraic
stacks, then we can construct the ``moduli stack'', that captures all
the information that we would like in a fine moduli space.

The idea of enlarging the category of algebraic varieties to study
moduli problems is not new. In fact A. Weil invented the concept of
abstract variety to give an algebraic construction of the Jacobian
of a curve.

These notes are an introduction to the theory of algebraic stacks.
I have tried to emphasize ideas and concepts through examples
instead of detailed proofs (I give references where these can be
found). In particular, section
\ref{sectionversus} is a detailed comparison between the moduli
\textit{scheme} and the moduli \textit{stack} of vector bundles.

First I will give a quick
introduction in subsection \ref{quick}, just to give some motivations
and get a flavour of the theory of algebraic stacks. 

Section \ref{sectionstacks} has a more detailed exposition.
There are mainly two ways of introducing stacks. We can
think of them as 2-functors (I learnt this approach from N. Nitsure
and C. Sorger, cf. subsection \ref{subsfunctors}), or as categories 
fibered on groupoids (This is the
approach used in the references, cf. subsection 
\ref{subsgroupoids}). From the first point of view
it is easier to see in which
sense stacks are generalizations of schemes, and the
definition looks more natural, so conceptually it seems more
satisfactory. But since the references use categories fibered on 
groupoids, after we present both points of view, we will
mainly use the second.

The concept of stack is merely a categorical concept. To do geometry
we have to add some conditions, and then we get the concept of
algebraic stack. This is done in subsection \ref{subsalgebraic}.

In subsection \ref{subsgroupspaces} we introduce a third point of view to
understand stacks: as groupoid spaces.

In subsection \ref{subsproperties} we define for algebraic stacks many of the
geometric properties that are defined for schemes (smoothness,
irreducibility, separatedness, properness, etc...). In
subsection \ref{subspoints} we introduce the concept of point and dimension of 
an algebraic stacks, and in subsection \ref{subssheaves} we define sheaves on 
algebraic stacks.

In section \ref{sectionversus} we study in detail the example 
of the moduli of
vector bundles on a scheme $X$, comparing the moduli stack with
the moduli scheme.

Appendix A is a brief introduction to Grothendieck topologies, sheaves
and algebraic spaces. In appendix B we define some notions related to
the theory of 2-categories.

\subsection{Quick introduction to algebraic stacks}
\label{quick}

We will start with an example: vector bundles (with fixed prescribed
Chern classes and rank) on a projective scheme $X$ over an
algebraically closed field $k$. What is the moduli stack
$\BBund$ of vector bundles on $X$?. I don't know a short answer to
this, but instead it is easy to define what is a morphism from a
scheme $B$ to the moduli stack $\BBund$. It is just a family of vector
bundles parametrized by $B$. More precisely, it is a vector bundle $V$
on $B\times X$, flat over $B$, such that the restriction to the slices
$b\times X$ have prescribed Chern classes and rank. In other words,
$\BBund$ has the property that we expect from a fine moduli space:
the set of morphisms $\Hom(B,\BBund)$ is equal to the set of families
parametrized by $B$.

We will say that a diagram
\begin{eqnarray}
\label{commdiag}
\xymatrix{
{B} \ar[r]^f \ar[rd]_{g} & {B'} \ar[d]^{g'} \\
 & {\BBund}
}
\end{eqnarray}
is commutative if the vector bundle $V$ on $B\times X$ corresponding
to $g$ is isomorphic to the vector bundle $(f\times \id_X)^*V'$, where
$V'$ is the vector bundle corresponding to $g'$.
Note that in general, if $L$ is a line bundle on $B$, then $V$ and
$V\otimes p^*_B L$ won't be isomorphic, and then the corresponding
morphisms from $B$ to $\BBund$ will be different, as opposed to what
happens with moduli schemes.

A $k$-point in the  stack $\BBund$ is a morphism
$u:\Spec k \to\BBund$, in other words, it is a vector
bundle $V$ on $X$, and we say that two points are isomorphic if they
correspond to isomorphic vector bundles. But we shouldn't think of $\BBund$
just as a set of points, it should be thought of as a 
category. The objects of $\BBund$ are points\footnote{
To be precise, we should consider also $B$-valued points,
for any scheme $B$, but we will only consider $k$-valued points for
the moment}, i.e. vector
bundles on $X$, and a morphism in $\BBund$ is an isomorphism of vector
bundles. This is the main difference between a scheme and an algebraic
stack: a scheme is a \textit{set} of points,
but in an algebraic stack is a \textit{category}, in fact a
\textit{groupoid} (i.e. a category in which all morphisms are isomorphisms). 
Each point comes with a group of
automorphisms. Roughly speaking, a scheme (or more generally, an
algebraic space \cite{Ar1}, \cite{K}) can be thought of as an
algebraic stack in which these groups of automorphisms are all trivial. 
If $p$ is the $k$-point in $\BBund$ corresponding to
a vector bundle $V$ on $X$, then the group of automorphisms
associated to $p$ is the group of vector bundle automorphisms of $V$.
This is why algebraic stacks are well suited to serve as moduli of
objects that have automorphisms.

An algebraic stack has an atlas. This is a scheme $U$ and a surjective
morphism $u:U \to \BBund$ (with some other properties). 
As we have seen, such a morphism $u$ is equivalent to a family of
vector bundles parametrized by $U$, and we say that $u$ is surjective
if for every vector bundle $V$ over $X$ there is at least one point in
$U$ whose corresponding vector bundle is isomorphic to $V$.
The existence of an atlas for an algebraic stack is the
analogue of the fact that for a scheme $B$ there is always an
\textit{affine} scheme $U$ and a surjective morphism $U \to B$ (if
$\{U_i\to B\}$ is a covering of $B$ by affine subschemes, take $U$ to
be the disjoint union $\coprod U_i$). Many local properties (smooth,
normal, reduced...) can be studied by looking at the atlas $U$. It is
true that in some sense an algebraic stack looks, locally, like a
scheme, but we shouldn't take this too far. For instance
the atlas of the classifying stack $BG$ (parametrizing principal
$G$-bundles, cf. example \ref{quotient}) is just a single point.
The dimension of an algebraic stack $\BBund$ will be defined as the
dimension of $U$ minus the relative dimension of the morphism $u$.
The dimension of an algebraic stack can be negative (for instance,
$\dim (BG)=-\dim(G)$).

A coherent sheaf $L$ on an algebraic stack $\BBund$ is a law that, 
for each morphism $g:B
\to \BBund$, gives a coherent sheaf $L_B$ on $B$, and for each
commutative diagram like (\ref{commdiag}), gives an isomorphism
between $f^* L_{B'}$ and $L_B$. The coherent sheaf $L_B$ should be
thought of as the pullback ``$g^*L$'' of $L$ under $g$ (the compatibility
condition for commutative diagrams is just the condition that
$(g'\circ f)^*L$ should be isomorphic to $f^* {g'}^* L$).

Let's look at another example: the moduli quotient (example 
\ref{quotient}). Let $G$ be an
affine algebraic group acting on $X$. For simplicity, assume that there is a 
normal subgroup $H$ of $G$ that acts trivially on $X$, and that
$\overline G=G/H$
is an affine group acting freely on $X$ and furthermore there is a 
quotient by this action 
$X \to B$ and this quotient is a principal $\overline G$-bundle.
We call $B=X/G$ the \textit{quotient scheme}. Each point 
corresponds to a
$G$-orbit of the action. But note that $B$ is also equal to the 
quotient $X/\overline G$, because $H$ acts trivially and then $G$-orbits are
the same thing as $\overline G$-orbits. We can say that the quotient scheme
``forgets'' $H$.

One can also define the \textit{quotient stack} $[X/G]$. Roughly
speaking, a point $p$ of $[X/G]$ again corresponds to a $G$-orbit of
the action, but now each
point comes with an automorphism group: given a point $p$ in $[X/G]$, 
choose a  point $x\in X$ in the orbit corresponding to $p$. The
automorphism group attached to $p$ is the stabilizer $G_x$ of $x$. With
the assumptions that we have made on the action of $G$, the
automorphism group of any point is always $H$.
Then the quotient stack $[X/G]$ is not a scheme, since the automorphism
groups are not trivial. The action of $H$ is trivial, but the
moduli stack still ``remembers'' that there was an action by $H$.
Observe that the stack $[X/\overline G]$ is not isomorphic to the stack
$[X/G]$ (as opposed to what happens with the quotient schemes). Since
the action of $\overline G$ is free on $X$, the automorphism group 
corresponding to
each point of $[X/\overline G]$ is trivial, and it can be shown 
that, with the assumptions that we made, $[X/\overline G]$ is
represented by the scheme $B$ (this terminology will be made precise
in section \ref{sectionstacks}).

\section{Stacks}
\label{sectionstacks}

\subsection{Stacks as 2-functors. Sheaves of sets.}
\label{subsfunctors}

Given a scheme $M$ over a base scheme $S$, we define its
(contravariant) functor of
points $\Hom_S(-,M)$
$$
\begin{array}{rccc}
{\Hom_S(-,M):} &{(\Sch /S)}& \longrightarrow &{(\Sets)} \\
               & {B} &\longmapsto &{\Hom_S(B,M)} 
\end{array}
$$
where $(\Sch /S)$ is the category of $S$-schemes, $B$ is an
$S$-scheme, and $\Hom_S(B,M)$ is the set of $S$-scheme morphisms.
If we give $(\Sch/S)$ the \'etale topology, $\Hom_S(-,M)$ is a sheaf.
A sheaf of sets on $(\Sch/S)$ with the \'etale topology is called a space.

Then schemes can be thought of as sheaves of sets. Moduli problems can
usually be
described by functors. We say that a sheaf of
sets $F$ is representable by a scheme $M$ if $F$ is isomorphic to the
functor of points $\Hom_S(-,M)$. The scheme $M$ is then called the
fine moduli scheme. 
Roughly speaking, this means that there is a one to one correspondence
between families of objects parametrized by a scheme $B$ and morphisms
from $B$ to $M$.

\begin{example}[Vector bundles]
\label{defvectorbundle}
\textup{ 
Let $X$ be a projective scheme over a Noetherian base $S$. 
We define the moduli functor 
$\FBund'$ of vector
bundles of fixed rank $r$ and Chern classes $c_i$ by sending the 
scheme $B$ to the set $\FBund'(B)$ of isomorphism classes of vector
bundles on $X\times B$, flat over $B$ with rank $r$ and whose 
restriction to the
slices $X\times \{b\}$ have Chern classes $c_i$. These vector bundles
should be thought of as families of vector bundles parametrized by $B$.
A morphism $f:B'\to B$ is sent to $\FBund'(f)=f^*:\FBund'(B) \to
\FBund'(B')$, the map of sets induced by the pullback. Usually we will
also fix a polarization $H$ in $X$ and restrict our attention to stable or
semistable vector bundles with respect to this polarization, and
then we consider the corresponding functors $\FBund^{\prime s}$ and
$\FBund^{\prime ss}$. 
}
\end{example}

\begin{example}[Curves]
\textup{
The moduli functor $M_g$ of smooth curves of genus $g$ over $S$ is the
functor that sends each scheme $B$ to the set $M_g(B)$ of isomorphism
classes of smooth and proper morphisms $C \to B$ (where $C$ is an 
$S$-scheme) whose fibers are geometrically connected curves of genus
$g$. Each morphism $f:B'\to B$ is sent to the map of sets induced by
the pullback $f^*$.
}
\end{example}

None of these examples are sheaves (then none of these are
representable), because of the presence
of automorphisms. They are just presheaves (=functors).
For instance, given a curve $C$ over $S$ with nontrivial
automorphisms, it is possible to construct a family $f:\SC \to B$ such
that every fiber of $f$ is isomorphic to $C$, but $\SC$ is not
isomorphic to $B \times C$. This implies that $M_g$ doesn't satisfy
the monopresheaf axiom.

This can be solved by taking the sheaf associated to the presheaf
(sheafification). In the examples, this amounts to change
isomorphism classes of families to equivalence classes of families,
when two families are equivalent if they are locally (using the
\'etale topology over the
parametrizing scheme $B$) isomorphic. In the case of vector bundles,
this is the reason why one usually declares two vector bundles $V$ and
$V'$ on $X \times B$ equivalent if $V\isom V'\otimes p_B^* L$ for some
line bundle $L$ on $B$. The functor obtained with this equivalence
relation is denoted $\FBund$ (and analogously for $\FBund^{s}$
and $\FBund^{ss}$).

Note that if two families $V$ and $V'$ are equivalent in this sense,
then they are locally isomorphic. The converse is only true if the
vector bundles are simple (only automorphisms are 
scalar multiplications). This will happen, for instance, if 
we are considering the functor $\FBund^{\prime s}$ of stable vector
bundles, since stable vector bundles are simple.
In general, if we want the functor to be a sheaf, we have to use
a weaker notion of equivalence, but this is not done because for 
other reasons there is only hope of obtaining a fine moduli space if
we restrict our attention to stable vector bundles.

Once this modification is made, there are some situations in which
these examples are representable (for instance, stable vector bundles 
on curves with coprime
rank and degree), but in general they 
will still not be representable,
because in general we don't have a universal family:

\begin{definition}[Universal family]
Let $F$ be a representable functor, and
let $\phi:F \to \Hom_S(-,X)$ be the isomorphism. The object of $F(X)$
corresponding to the element $\id_X$ of $\Hom_S(X,X)$ is called the
universal family.
\end{definition}
  
\begin{example}[Vector bundles]
\textup{
If $V$ is a universal vector bundle (over $S\times M$, where $M$
is the fine moduli space), it has the property that for 
any family $W$
of vector bundles (i.e. $W$ is a vector bundle over $X\times B$ for
some parameter scheme $B$) there exists a morphism $f:B\to M$ such that 
$(f\times \id_X)^* V$ is equivalent to $W$.}
\end{example}

When a moduli functor $F$ is not representable and then there is no scheme
$X$ whose functor of points is isomorphic to $F$, one can still try to
find a scheme $X$ whose functor of points is an approximation to $F$
in some sense. There are two different notions:

\begin{definition}[Corepresents]
\textup{\cite[p. 60]{S}, \cite[def 2.2.1]{HL}}.
We say that a scheme $M$ corepresents the functor $F$ if there 
is a natural
transformation of functors $\phi:F \to \Hom_S(-,M)$ such that
\begin{itemize}

\item Given another scheme $N$ and a natural transformation $\psi:F
\to \Hom_S(-,N)$, there is a unique natural transformation $\eta:
\Hom_S(-,M)\to \Hom_S(-,N)$ with $\psi= \eta \circ \phi$.
$$
\xymatrix{
{F} \ar[d]^{\phi} \ar[rd]^{\psi} \\
 {\Hom_S(-,M)} \ar[r]^{\eta} &  \Hom_S(-,N)\\ 
}
$$
\end{itemize}
\end{definition}

This characterizes $M$ up to unique isomorphism. Let $(\Sch/S)'$ be the
functor category, whose objects are contravariant functors from
$(\Sch/S)$ to $(Sets)$ and whose morphisms are natural transformation
of functors. Then  $M$ represents $F$ iff $\Hom_S(Y,M)=
\Hom_{(\Sch/S)'}(\SY,F)$ for all schemes $Y$, where
$\SY$ is the functor represented by $Y$. On the other hand,
one can check that $M$ corepresents $F$ iff $\Hom_S(M,Y)=
\Hom_{(\Sch/S)'}(F,\SY)$ for all schemes $Y$.
If $M$ represents $F$, then it corepresents it, but the converse is
not true. From now on we will usually denote a scheme and the functor 
that it represents by the same letter.

\begin{definition}[Coarse moduli]
A scheme $M$ is called a coarse moduli scheme if it corepresents $F$
and furthermore

\begin{itemize}
\item For any algebraically closed field $k$, the map 
$\phi(k):F(\Spec k) \to \Hom_S(\Spec k, M)$ is bijective.

\end{itemize}
\end{definition}

In both cases, given a family of objects parametrized by $B$ we get a 
morphism
from $B$ to $M$, but we don't require the converse to be true.

\begin{example}[Vector bundles]
\label{vb1}
\textup{
There is a scheme $\Bund^{ss}$ that corepresents $\FBund^{ss}$. It 
fails to be a
coarse moduli scheme because its closed points are in one to one
correspondence with S-equivalence classes of vector bundles, and not
with isomorphisms classes of vector bundles. Of course, this can be
solved `by hand' by modifying the functor and considering two vector
bundles equivalent if they are S-equivalent. Once this modification 
is done, $\Bund^{ss}$ is
a coarse moduli space.
}

\textup{
But in general $\Bund^{ss}$ doesn't represent the moduli functor
$\FBund^{ss}$. The reason for this is that vector bundles have always
nontrivial automorphisms (multiplication by scalar), but the moduli
functor doesn't record information about automorphisms: recall that to
a scheme $B$ it associates just the set of equivalence classes of vector
bundles. To record the automorphisms of these vector bundles, we define
$$
\begin{array}{rccc}
\BBund: & (\Sch/S) & \longrightarrow & (\groupoids) \\
        &    B     & \longmapsto     & \BBund(B)
\end{array}
$$
where $\BBund(B)$ is the category whose objects are vector bundles $V$
on $X\times B$ of rank $r$ and with fixed Chern classes (note that the
objects are vector bundles, not isomorphism classes of vector
bundles), and whose 
morphisms are vector bundle isomorphisms (note that we use
isomorphisms of vector bundles, not S-equivalence nor equivalence
classes as before).
This defines a 2-functor between the 2-category associated to
$(\Sch/S)$ and the 2-category $(\groupoids)$ .}
\end{example}

\begin{definition}
Let $(\groupoids)$ be the 2-category whose objects are
groupoids, 1-morphisms are functors between groupoids, and 2-morphisms
are natural transformation between these functors.
A presheaf in groupoids (also called a quasi-functor) is a
contravariant 2-functor $\SF$ from $(\Sch/S)$ to $(\groupoids)$. 
For
each scheme $B$ we have a groupoid $\SF(B)$ and for each morphism $f:B'\to
B$ we have a natural transformation of functors $\SF(f)$ that is
denoted by $f^*$ (usually it is actually defined by a pullback).
\end{definition}

\begin{example}[Vector bundles]
\label{bbund}
\textup{\cite[1.3.4]{La}.
$\BBund$ is a presheaf. For each object $B$ of 
$(\Sch/S)$ it gives the
groupoid $\BBund(B)$ that we have defined in example \ref{vb1}. 
For each 1-morphism
$f:B' \to B$ it gives the functor $F(f)=f^*:\BBund(B)\to \BBund(B')$
given by pull-back, and for every diagram
\begin{eqnarray}
\label{compo}
B'' \stackrel{g}\longrightarrow B' \stackrel{f}\longrightarrow B
\end{eqnarray}
it gives a natural transformation of functors (a 2-isomorphism)
$\epsilon_{g,f}:g^*\circ f^* \to (f\circ g)^*$. This is the only
subtle point. First recall that the pullback $f^*V$ of a vector bundle
(or more generally, any fiber product) is not uniquely defined: it is
only defined up to unique isomorphism. First choose once and for all a
pullback $f^*V$ for each $f$ and $V$. Then, given a diagram like
\ref{compo}, in principle $g^*(f^*V)$ and $(f\circ g)^*V$ are not the
same, but (because both solve the same universal problem) there is a
canonical isomorphism (the unique isomorphism of the universal
problem) $g^*(f^*V) \to (f\circ g)^*V$ between them,
and this defines the natural transformation of functors 
$\epsilon_{g,f}:g^*\circ f^* \to (f\circ g)^*$. By a slight abuse of
language, usually we won't write explicitly these isomorphisms
$\epsilon_{g,f}$, and we will write $g^*\circ f^* = (f\circ g)^*$. 
Since they are uniquely defined this will cause no
ambiguity.}
\end{example}

Now we will define the concept of stack. First we have to choose a
Grothendieck topology on $(Sch/S)$, either the \'etale or the fppf
topology. Later on, when we define algebraic stack, the \'etale 
topology will lead to the definition of a
Deligne-Mumford stack (\cite{DM}, \cite{Vi}, \cite{E}), and the fppf to
an Artin stack (\cite{La}). For the moment we will give a unified
description.

In the following definition, to simplify notation we denote by $X|_i$
the pullback $f^*_i X$ where $f_i:U_i \to U$ and $X$ is an object of 
$\SF(U)$, and by $X_i|_{ij}$ the
pullback $f^*_{ij,i} X_i$ where $f_{ij,i}:U_i \times_U U_j \to U_i$
and $X_i$ is an object of $\SF(U_i)$. We will also use the obvious
variations of this convention, and will simplify the notation using
remark \ref{B2}.

\begin{definition}[Stack]
\label{sheaf}
A stack is a sheaf of groupoids, i.e. a 2-functor (presheaf) that
satisfies the following sheaf axioms. Let $\{U_i \to U\}_{i\in I}$ be
a covering of $U$ in the site $(\Sch/S)$. Then
\begin{enumerate}
\item (Glueing of morphisms) If $X$ and $Y$ are two objects of 
$\SF(U)$, and $\varphi_i:X|_i\to
Y|_i$ are morphisms such that $\varphi_i|_{ij}=\varphi_j|_{ij}$, then
there exists a morphism $\eta:X\to Y$ such that $\eta|_i=\varphi_i$.

\item (Monopresheaf) If $X$ and $Y$ are two objects of $\SF(U)$, 
and $\varphi:X\to
Y$, $\psi:X \to Y$ are morphisms such that $\varphi|_i=\psi|_i$, then
$\varphi = \psi$.

\item \label{sheafthree} (Glueing of objects) If $X_i$ are objects 
of $\SF(U_i)$ and 
$\varphi_{ij}:X_j|_{ij}
\to X_i|_{ij}$ are morphisms satisfying the cocycle condition 
$\varphi_{ij}|_{ijk}\circ \varphi_{jk}|_{ijk}= \varphi_{ik}|_{ijk}$,
then there exists an object $X$ of $\SF(U)$ and $\varphi_i:X|_i
\stackrel{\isom}\to X_i$ such that $\varphi_{ji}\circ \varphi_i|_{ij}=
\varphi_j|_{ij}$.
\end{enumerate}
\end{definition}

Let's stop for a moment and look at how we have enlarged the category
of schemes by defining the category of stacks. We can draw the following
diagram
$$
\xymatrix{
& {Algebraic\,Stacks} \ar[r] & {Stacks} \ar[r] &{Presheaves\,of\,
groupoids} \\
{\Sch/S} \ar[r] \ar[ur] &
{Algebraic\,Spaces} \ar[r] \ar[u]
&{Spaces} \ar[r] \ar[u] &{Presheaves\,of\,sets} \ar[u]
}
$$
where $A \to B$ means that the category $A$ is a subcategory $B$. 
Recall that a presheaf of sets is just
a functor from $(\Sch/S)$ to the category $(\Sets)$, a presheaf of 
groupoids is
just a 2-functor to the 2-category $(\groupoids)$. A sheaf (for example an
space or a stack) is a presheaf that satisfies the sheaf axioms
(these axioms are slightly different in the context of categories or
2-categories), and if this sheaf satisfies some geometric conditions
(that we haven't yet specified), we will have an algebraic stack or
algebraic space.

\subsection{Stacks as categories. Groupoids}
\label{subsgroupoids}

There is an alternative way of defining a stack. From this point of
view a stack will be a category, instead of a functor.

\begin{definition}
A category over $(\Sch/S)$ is a category $\SF$ and a covariant functor
$p^{}_\SF:\SF \to (\Sch/S)$. If $X$ is an object (resp. $\phi$ is a
morphism) of $\SF$, and $p^{}_\SF(X)=B$ (resp. $p^{}_\SF(\phi)=f$), then we
say that $X$ lies over $B$ (resp. $\phi$ lies over $f$).
\end{definition}

\begin{definition}[Groupoid]
A category $\SF$ over $(\Sch/S)$ is called a category fibered on
groupoids (or just groupoid) if 
\begin{enumerate}
\item \label{groupoidone} For every $f:B'\to B$ in $(\Sch/S)$ and every object $X$ with 
$p^{}_\SF(X)=B$, there exists at least one object $X'$ and a morphism 
$\phi:X'\to
X$ such that $p^{}_\SF(X')=B'$ and $p^{}_\SF(\phi)=f$.
$$
\xymatrix{
{X'} \ar@{-->}[r]^{\phi} \ar@{-->}[d] & {X} \ar[d] \\
{B'} \ar[r]^{f} & {B} }
$$

\item \label{groupoidtwo} For every diagram
$$\xymatrix{
{X_3} \ar[rr]^{\psi} \ar[dd]& & {X_1} \ar[dd] \\
 & {X_2} \ar[ru]^{\phi} \ar[dd] \\
{B_3} \ar '[r][rr]^{f\circ f'} \ar[rd]_{f'} & & {B_1} \\
 & {B_2} \ar[ru]_f 
}
$$
(where $p^{}_\SF(X_i)=B_i$, $p^{}_\SF(\phi)=f$, $p^{}_\SF(\psi)=f\circ f'$),
there exists a unique $\varphi:X_3 \to X_2$ with $\psi=\phi\circ
\varphi$ and $p^{}_\SF(\varphi)=f'.$
\end{enumerate}
\end{definition}

Condition \ref{groupoidtwo} implies that the object $X'$ whose existence
is asserted in condition \ref{groupoidone} is unique up to canonical
isomorphism. For each $X$ and $f$ we choose once and for all such an
$X'$ and call it $f^*X$.
Another consequence of condition \ref{groupoidtwo} is 
that $\phi$ is an isomorphism
if and only if $p^{}_\SF(\phi)=f$ is an isomorphism.

Let $B$ be an object of $(\Sch/S)$. We define $\SF(B)$, the fiber of
$\SF$ over $B$, to be the subcategory of $\SF$ whose objects lie over
$B$ and whose morphisms lie over $\id_B$. It is a groupoid.

The association  $B\to \SF(B)$ in fact defines a presheaf of groupoids
(note that the 2-isomorphisms $\epsilon_{f,g}$ required in the
definition of presheaf of groupoids are well defined thanks to 
condition \ref{groupoidtwo}). Conversely, given a presheaf of 
groupoids $\SG$ on
$(Sch/S)$, we can define the category $\SF$ whose objects are pairs
$(B,X)$ where $B$ is an object of $(\Sch/S)$ and $X$ is an object of
$\SG(B)$, and whose morphisms $(B',X')\to (B,X)$ are pairs
$(f,\alpha)$ where $f:B'\to B$ is a morphism in $(Sch/S)$ and
$\alpha:f^* X \to X'$ is an isomorphism, where $f^*=\SG(f)$.
This gives the relationship between both points of view.

\begin{example}[Stable curves]
\label{defstablecurve}
\textup{\cite[def 1.1]{DM}. 
Let $B$ be an $S$-scheme. Let $g\geq 2$. A stable curve of genus $g$
over $B$ is a proper and flat morphism $\pi:C \to B$ whose geometric
fibers are reduced, connected and one-dimensional schemes $C_b$ such
that 
\begin{enumerate}
\item The only singularities of $C_b$ are ordinary double points.
\item If $E$ is a non-singular rational component of $C_b$, then $E$
meets the other components of $C_b$ in at least 3 points.
\item $\dim H^1(\SO_{C_b})=g$.
\end{enumerate}
Condition 2 is imposed so that the automorphism group of $C_b$ is finite.
A stable curve over $B$ should be thought of as a family of stable
curves (over $S$) parametrized by $B$.}

\textup{
We define $\overline\SM_g$, the groupoid over $S$ whose objects
are stable curves over $B$ and whose morphisms are Cartesian diagrams
$$
\xymatrix{
{X'} \ar[r] \ar[d] & {X} \ar[d] \\
{B'} \ar[r] & {B}}
$$}
\end{example}

\begin{example}[Quotient by group action]
\label{quotient}
\textup{\cite[1.3.2]{La}, \cite[example 4.8]{DM}, 
\cite[example 2.2]{E}.
Let $X$ be an $S$-scheme (assume all schemes are Noetherian), 
and $G$ an affine flat group $S$-scheme acting on
the right on $X$. We define the groupoid $[X/G]$ whose objects are
principal $G$-bundles  $\pi:E\to B$ together with a $G$-equivariant
morphism $f:E\to X$. A morphism is Cartesian diagram
$$
\xymatrix{
{E'} \ar[r]^{p} \ar[d]_{\pi'} & {E} \ar[d]_{\pi} \\
{B'} \ar[r] & {B}}
$$
such that $f\circ p= f'$.}

\end{example}

\begin{definition}[Stack]
A stack is a groupoid that satisfies
\begin{enumerate}
\item (\textit{Prestack}). For all scheme $B$ and pair of objects $X$, 
$Y$ of $\SF$ over $B$,
the contravariant functor
$$
\begin{array}{rccc}
\Iso_B(X,Y): & (\Sch/B)& \longrightarrow & (\Sets) \\
 & (f:B'\to B) & \longmapsto & \Hom(f^*X,f^*Y) 
\end{array}
$$
is a sheaf on the site $(Sch/B).$

\item Descent data is effective (this is just condition
\ref{sheafthree} in the
definition \ref{sheaf} of sheaf).
\end{enumerate}
\end{definition}

\begin{example}
\textup{
If $G$ is smooth and affine, the groupoid $[X/G]$ is a stack 
\cite[2.4.2]{La}, \cite[example 7.17]{Vi}, \cite[prop 2.2]{E}.
Then also $\overline\SM_g$ (cf. example \ref{defstablecurve}) 
is a stack, because it is isomorphic to a
quotient stack of a subscheme of a Hilbert scheme by $\pgl$ 
\cite[thm 3.2]{E}, \cite{DM}. 
The groupoid $\BBund$ defined in example \ref{defvectorbundle} 
is also a stack \cite[2.4.4]{La}.}
\end{example}
From now on we will mainly use this approach. Now we will give some
definitions for stacks.
\medskip

\textbf{Morphisms of stacks}. A morphism of stacks $f:\SF \to \SG$ 
is a functor between the
categories, such that $p_\SG \circ f= p^{}_\SF$.
A commutative diagram of stacks is a diagram
$$
\xymatrix{
 & {\SG} \ar[rd]^g \ar@2[d]^{\alpha} \\
{\SF} \ar[ur]^f \ar[rr]_h & &{\SH}
}
$$
such that $\alpha:g\circ f \to h$ is an isomorphism of functors.
If $f$ is an equivalence of categories, then we say that the stacks 
$\SF$ and $\SG$ are isomorphic.
We denote by $\Hom_S(\SF,\SG)$ the category whose objects are morphisms
of stacks and whose morphisms are natural transformations.

\medskip
\textbf{Stack associated to a scheme}. Given a scheme $U$ over $S$, 
consider the category $(\Sch/U)$. Define
the functor $p^{}_U:(\Sch/U)\to (\Sch/S)$ which sends the $U$-scheme
$f:B\to U$ to the composition $B\stackrel{f}\to U \to S$. Then
$(\Sch/U)$ becomes a stack. Usually we denote this stack also by $U$.
From the point of view of 2-functors, the stack associated to $U$ is
the 2-functor that for each scheme $B$ gives the category whose
objects are the elements of the set $\Hom_S(B,U)$, and whose only
morphisms are identities. 

We say that a stack is represented by a scheme $U$ when it is
isomorphic to the stack associated to $U$. We have the following
very useful lemmas:

\begin{lemma}
\label{nonrepresentable}
If a stack has an
object with an automorphism other that the identity, then the stack  
cannot be represented by a scheme. 
\end{lemma}

\begin{proof}
In the definition of stack associated with a scheme we see that 
the only automorphisms are identities.
\end{proof}

\begin{lemma}
\label{yoneda}
\textup{\cite[7.10]{Vi}}.
Let $\SF$ be a stack and $U$ a scheme. The functor
$$
u:\Hom_S(U,\SF) \to \SF(U)
$$
that sends a morphism of stacks $f:(\Sch/U)\to \SF$ to $f(\id_U)$ is
an equivalence of categories. 
\end{lemma}

\begin{proof}
Follows from Yoneda lemma
\end{proof}

This useful
observation that we will use very often means that an object of $\SF$
that lies over $U$ is equivalent to a morphism (of stacks) from $U$ to
$\SF$.

\medskip
\textbf{Fiber product}. Given two morphisms $f_1:\SF_1\to \SG$,
$f_2:\SF_2\to \SG$, we define a new stack $\SF_1 \times_\SG \SF_2$
(with projections to $\SF_1$ and $\SF_2$) as follows.
The objects are triples $(X_1,X_2,\alpha)$ where $X_1$ and $X_2$ are
objects of $\SF_1$ and $\SF_2$ that lie over the same scheme $U$, and
$\alpha: f_1(X_1)\to f_2(X_2)$ is an isomorphism in $\SG$
(equivalently,  $p_\SG(\alpha)=\id_U$).
A morphism from $(X_1,X_2,\alpha)$ to $(Y_1,Y_2,\beta)$ is a pair
$(\phi_1,\phi_2)$ of morphisms $\phi_i:X_i\to Y_i$ that lie over the
same morphism of schemes $f:U \to V$, and such that $\beta \circ
f_1(\phi_1) = f_2(\phi_2)\circ \alpha$.
The fiber product satisfies the usual universal property.

\medskip
\textbf{Representability}. A stack $\SX$ is said to be representable by an
algebraic space (resp. scheme) if there is an algebraic space
(resp. scheme) $X$ such that the stack associated to $X$ is isomorphic
to $\SX$.
If ``P'' is a property of algebraic spaces (resp. schemes) and $\SX$ 
is a representable stack, we will say that $\SX$ has ``P'' iff $X$
has ``P''.

A morphism of stacks $f:\SF\to \SG$ is said to be representable if for
all objects $U$ in $(\Sch/S)$ and morphisms $U\to \SG$, the fiber
product stack $U\times_\SG \SF$ is representable by an algebraic
space.
Let ``P'' is a property of morphisms of schemes that is local in nature
on the target for the topology chosen on $(\Sch/S)$ (\'etale or
fppf), and it is stable under arbitrary base change. For instance:
separated, quasi-compact, unramified, flat, smooth, \'etale, surjective,
finite type, locally of finite type,... Then we say that
$f$ has ``P'' if for every $U\to \SG$, the pullback $U\times_\SG \SF
\to U$ has ``P'' (\cite[p.17]{La}, \cite[p.98]{DM}).

\medskip
\textbf{Diagonal}. Let $\Delta_\SF:\SF \to \SF\times_S \SF$ be the
obvious diagonal morphism. A morphism from a scheme $U$ to $\SF 
\times_S \SF$ is equivalent to two objects $X_1$, $X_2$ of
$\SF(U)$. Taking the fiber product of these we have
$$
\xymatrix{
{\Iso_U(X_1,X_2)} \ar[r] \ar[d]& {\SF} \ar[d]^{\Delta_\SF} \\
{U} \ar[r]^{(X_1,X_2)} & {\SF\times_S \SF}}
$$
hence the group of automorphisms of an object is encoded in the
diagonal morphism.

\begin{proposition}
\label{diag}
\textup{\cite[cor 2.12]{La}, \cite[prop 7.13]{Vi}}.
The following are equivalent
\begin{enumerate}
\item The morphism $\Delta_\SF$ is representable.

\item The stack $\Iso_U(X_1,X_2)$ is representable for all $U$, $X_1$
and $X_2$.

\item For all scheme $U$, every morphism $U\to \SF$ is representable.

\item For all schemes $U$, $V$ and morphisms $U\to \SF$ and $V\to \SF$,
the fiber product $U\times_\SF V$ is representable.
\end{enumerate}
\end{proposition}

\begin{proof}

The implications $1 \Leftrightarrow 2$ and $3 \Leftrightarrow 4$
follow easily from the definitions.

$1 \Rightarrow 4$) Assume that $\Delta_\SF$ is representable. We have
to show that $U\times_\SF V$ is representable for any $f:U\to \SF$ and
$g:V\to \SF$. Check that the following diagram is Cartesian
$$
\xymatrix{
{U\times_\SF V} \ar[r] \ar[d]& {\SF}\ar[d]^{\Delta_\SF}\\
U\times_S V \ar[r]^{f\times g} &{\SF\times_S \SF}}
$$
Then $U\times_\SF V$ is representable.

$1 \Leftarrow 4$) First note that the Cartesian diagram defined by
$h:U\to \SF\times_S \SF$ and $\Delta_\SF$ factors as follows
$$
\xymatrix{
{U\times^{}_{\SF\times_S \SF} \SF} \ar[r] \ar[d] & 
{U\times^{}_\SF U} \ar[r] \ar[d] &{\SF} \ar[d] \\
{U} \ar[r]^{\Delta_U} & {U\times_S U} \ar[r] &
{\SF\times_S \SF}}
$$
Both squares are Cartesian and by hypothesis $U\times_\SF U$ is
representable, then $U\times^{}_{\SF\times_S \SF} \SF$ is also
representable.

\end{proof}

\subsection{Algebraic stacks}
\label{subsalgebraic}

Now we will define the notion of algebraic stack. As we have said,
first we have to choose a topology on $(\Sch/S)$. Depending of whether
we choose the \'etale or fppf topology, we get different notions.

\begin{definition}[Deligne-Mumford stack] 
Let $(\Sch/S)$ be the category of
$S$-schemes with  the \'etale topology. Let $\SF$ be a stack. Assume
\begin{enumerate}
\item The diagonal $\Delta_\SF$ is representable, quasi-compact and
separated. 

\item There exists a scheme $U$ (called atlas) and an \'etale
surjective morphism
$u:U\to \SF$.
\end{enumerate}
Then we say that $\SF$ is a Deligne-Mumford stack.
\end{definition}

The morphism of stacks $u$ is representable because of proposition
\ref{diag} and the fact that the diagonal $\Delta_\SF$ is
representable. Then the notion of \'etale is well defined for $u$.
In \cite{DM} this was called an algebraic stack. In the literature,
algebraic stack usually refers to Artin stack (that we will define
later). To avoid confusion, we will use ``algebraic stack'' only when
we refer in general to both notions, and we will use
``Deligne-Mumford'' or ``Artin'' stack when we want to be specific.

Note that the definition of Deligne-Mumford stack is the same as the
definition of algebraic space, but in the context of stacks instead of
spaces.
As with schemes a stack such that the diagonal $\Delta_\SF$ is 
quasi-compact and
separated is called quasi-separable. We always assume this technical
condition, as it is usually done both with schemes and algebraic
spaces.

Sometimes it is difficult to find explicitly an \'etale atlas, and the
following proposition is useful.

\begin{proposition}
\label{represen}
\textup{\cite[thm 4.21]{DM}, \cite{E}}.
Let $\SF$ be a stack over the \'etale site $(\Sch/S)$. Assume
\begin{enumerate}
\item The diagonal $\Delta_\SF$ is representable, quasi-compact,
separated and \textbf{unramified}.

\item There exists a scheme $U$ of finite type over $S$ and a
\text{smooth} surjective morphism
$u:U\to \SF$.
\end{enumerate}
Then $\SF$ is a Deligne-Mumford stack.
\end{proposition}

Now we define the analogue for the fppf topology \cite{Ar2}.

\begin{definition}[Artin stack]
Let $(\Sch/S)$ be the category of
$S$-schemes with  the fppf topology. Let $\SF$ be a stack. Assume
\begin{enumerate}
\item The diagonal $\Delta_\SF$ is representable, quasi-compact and
separated. 

\item There exists a scheme $U$ (called atlas) and a smooth 
(hence locally of finite type) and 
surjective morphism
$u:U\to \SF$.
\end{enumerate}
Then we say that $\SF$ is an Artin stack.
\end{definition}

For propositions analogous to proposition \ref{represen} see \cite[4]{La}.

\begin{proposition}
\textup{\cite[prop 7.15]{Vi}, \cite[lemme 3.3]{La}}.
If $\SF$ is a
Deligne-Mumford
(resp. Artin) stack, then the diagonal
$\Delta_\SF$ is unramified (resp. finite type).
\end{proposition}

Recall that $\Delta_\SF$ is unramified (resp. finite type) if for
every scheme $B$ and objects $X$, $Y$ of $\SF(B)$, the morphism
$\Iso_B(X,Y)\to U$ is unramified (resp. finite type). If $B=\Spec S$
and $X=Y$, then this means that the automorphism group of $X$ is
discrete and reduced for a Deligne-Mumford stack, and it just of
finite type for an Artin stack.

\begin{example}[Vector bundles]
\label{quotconstruction}
\textup{The stack $\BBund$ is an Artin stack, locally of finite type 
\cite[4.14.2.1]{La}. The atlas is constructed as follows. Let 
Let $P^H_{r,c_i}$ be the Hilbert polynomial corresponding to
sheaves on $X$ with rank $r$ and Chern classes $c_i$.
Let $\Quot(\SO(-m)^{\oplus N}, P^H_{r,c_i})$ be the Quot scheme
parametrizing quotients of sheaves on $X$
\begin{eqnarray}
\label{quotmap}
\SO(-m)^{\oplus N} \surj V,
\end{eqnarray}
where $V$ is a coherent sheaf on $X$ with Hilbert polynomial
$P^H_{r,c_i}$. Let $R_{N,m}$ be the subscheme corresponding to quotients
(\ref{quotmap}) such that $V$ is a vector
bundle with $H^p(V(m))=0$ for $p>0$ and the morphism (\ref{quotmap})
induces an isomorphism on global sections
$$
H^0(\SO)^{\oplus N} \stackrel{\isom}{\too} H^0(V(m)).
$$
The scheme $R^{}_{N,m}$ has a universal vector bundle, induced from the
universal bundle of the Quot scheme, and then there is a morphism
$u^{}_{N,m}: R^{}_{N,m}\to \BBund$. 
Since $H$ is ample, for every vector bundle $V$, 
there exist integers $N$ and $m$ such that
$R_{N,m}$ has a point whose corresponding quotient is $V$, and then if
we take the infinite disjoint union of these morphisms 
we get a surjective morphism
$$
u: \Big( \coprod_{N,m>0} R^{}_{N,m}\Big) \too \BBund. 
$$
It can be shown that this morphism is smooth, and then it gives an
atlas. Each scheme $R_{N,m}$ is of finite type, so the union is
locally of finite type, which in turn implies that the stack $\BBund$
is locally of finite type.
}

\end{example}

\begin{example}[Quotient by group action]
\label{atlasquotient}
\textup{The stack $[X/G]$ is an Artin stack
\cite[4.14.1.1]{La}. 
If $G$ is smooth, an atlas is defined as follows (for more general
$G$, see \cite[4.14.1.1]{La}):
Take the trivial principal $G$-bundle
$X\times G$ over $X$, and let the map $f:X\times G \to X$ be the
action of the group. This defines an object of $[X/G](X)$, and by
lemma \ref{yoneda}, it defines a morphism $u:X\to [X/G]$. It is
representable, because if $B$ is a scheme and $g:B\to [X/G]$ is the 
morphism corresponding
to a principal $G$-bundle $E$ over $B$ with an equivariant morphism
$f:E\to X$, then $B\times_{[X/G]}X$ is isomorphic to the scheme $E$,
and in fact we have a Cartesian diagram
$$
\xymatrix{
{E} \ar[r]^{f} \ar[d]_{\pi} & {X} \ar[d]_{u} \\
{B} \ar[r]^{g} & {[X/G].}
}
$$
The morphism $u$ is surjective and smooth because $\pi$ is surjective
and smooth for every $g$ (if $G$ is not smooth, but only separated,
flat and of finite presentation, then $u$ is not an atlas, but if we
apply the representation theorem \cite[thm 4.1]{La}, we conclude that
there is a smooth atlas).}

\textup{
If either $G$ is \'etale 
over $S$ (\cite[example 4.8]{DM}) or 
the stabilizers of the geometric points of $X$ are finite and reduced
(\cite[example 7.17]{Vi}),
then $[X/G]$ is a Deligne-Mumford stack. In particular
$\overline\SM_g$ is a Deligne-Mumford stack.}

\textup{Note that if the action is not free, then $[X/G]$ is not
representable by lemma \ref{nonrepresentable}. On the other hand, 
if there is a scheme $Y$ such that  $X \to Y$ is a principal 
$G$-bundle, then $[X/G]$ is represented by $Y$.} 

\textup{Let $G$ be a reductive group acting on $X$. Let $H$ be an
ample line bundle on $X$, and assume that the action is polarized.
Let $X^s$ and $X^{ss}$ be the subschemes of stable and semistable
points. Let $Y=X\gitq G$ be the GIT quotient. 
Recall that there is a 
good quotient $X^{ss}\to Y$, and that the restriction to the stable 
part $X^s\to Y$ is a principal bundle. There is a natural morphism 
$[X^{ss}/G] \to X^{ss}\gitq G$. By the previous remark, the restriction
$[X^s/G] \to Y^s$ is an isomorphism of stacks.}

\textup{If $X=S$ (with trivial action of $G$ on $S$), then $[S/G]$ is 
denoted $BG$, the classifying groupoid of principal $G$-bundles.}

\end{example}

\subsection{Algebraic stacks as groupoid spaces}
\label{subsgroupspaces}

We will introduce a third equivalent definition of stack.
First consider a category $C$. Let $U$ be the set of objects and $R$
the set of morphisms. The axioms of a category give us four maps of
sets
$$
\xymatrix{
{R} \ar@<0.5ex>[r]^{s} 
\ar@<-0.5ex>[r]_{t} & {U} \ar[r]^{e} & {R}}
\qquad
\xymatrix{
\save[]+<-5.5ex,-0.55ex>*{R\times^{}_{s,U,t} R}\restore \ar[r]^{m} & {R}}
$$
where $s$ and $t$ give the source and target for each morphism, $e$
gives the identity morphism, and $m$ is composition of morphisms.
If the category is a groupoid then we have a fifth morphism
$$
\xymatrix{{R} \ar[r]^i & {R}}
$$
that gives the inverse. These maps satisfy
\begin{enumerate}
\item $s\circ e= t\circ e = \id_R$, $s\circ i=t$, $t\circ i=s$,
$s\circ m=s\circ p_2$, $t\circ m=t\circ p_1$.

\item \textit{Associativity}. $m\circ (m\times \id_R)=m\circ 
(\id_R \times m)$.

\item \textit{Identity}. Both compositions
$$
R=R\times^{}_{s,U} U=U\times^{}_{U,t}R
\xymatrix{
{}\ar@<0.5ex>[r]^{\id_R \times e} 
\ar@<-0.5ex>[r]_{e \times \id_R} & {}}
R\times^{}_{s,U,t} R
\xymatrix{
{}\ar[r]^{m} & {R}}
$$
are equal to the identity map on $R$.

\item \textit{Inverse}. $m\circ (i\times \id_R)= e\circ s$, 
$m\circ (\id_R \times i)= e\circ t$.
\end{enumerate}

\begin{definition}[Groupoid space]
\textup{\cite[1.3.3]{La}, \cite[pp. 668--669]{DM}}.
A groupoid space is a pair of spaces (sheaves of sets) $U$, $R$, with
five morphisms $s$, $t$, $e$, $m$, $i$ with the same properties as
above.
\end{definition}

\begin{definition}
\textup{\cite[1.3.3]{La}}.
Given a groupoid space, define the groupoid over $(\Sch/S)$ as the
category $[R,U]'$ over $(\Sch/S)$ whose objects over the scheme $B$
are elements of the set $U(B)$ and whose morphisms over $B$ are
elements of the set $R(B)$. Given $f:B' \to B$ we define a functor
$f^*: [R,U]'(B) \to [R,U]'(B')$ using the maps $U(B) \to U(B')$ and
$R(B) \to R(B')$. 
\end{definition}

The groupoid $[R,U]'$ is in general only a prestack. We denote by
$[R,U]$ the associated stack. The stack $[R,U]$ can be thought of as 
the sheaf associated to the presheaf of groupoids $B \mapsto
[R,U]'(B)$ (\cite[2.4.3]{La}).

\begin{example}[Quotient by group action]
\textup{Let $X$ be a scheme and $G$ an affine group scheme. We denote
by the same letters the associated spaces (functors of points). We
take $U=X$ and $R=X\times G$. Using the group action we can define the
five morphisms ($t$ is the action of the group, $s=p_1$, $m$ is the
product in the group, $e$ is defined with the identity of $G$, and $i$
with the inverse).}

\textup{
The objects of $[X\times G,X]'(B)$ are morphisms $f:B\to
X$. Equivalently, they are trivial principal $G$-bundles $B\times G$
over $B$ and a map $B\times G \to X$ defined as the composition of the
action of $G$ and $f$. The stack $[X\times G,X]$ is isomorphic to $[X/G]$.}
\end{example}

\begin{example}[Algebraic stacks]
\textup{
Let $R$, $U$ be a groupoid space such that $R$ and $U$ are algebraic
spaces, locally of finite presentation (equivalently locally of finite
type if $S$ is noetherian). Assume that the morphisms $s$, $t$ are
flat, and that $\delta=(s,t):R\to U\times_S U$ is separated and
quasi-compact. Then $[R,U]$ is an Artin stack, locally of finite type 
(\cite[cor 4.7]{La}).}

\textup{
In fact, any Artin stack $\SF$ can be defined in this fashion. The
algebraic space $U$ will be the atlas of $\SF$, and we set 
$R=U\times_\SF U$. The morphisms $s$ and $t$ are the two projections,
$i$ exchanges the factors, $e$ is the diagonal, and $m$ is defined by
projection to the first and third factor.}
\end{example}

Let $\delta:R\to U\times_S  U$ be an equivalence relation in the
category of spaces. One can define a groupoid space, and $[R,U]$ is to
be thought of as the stack-theoretic quotient of this equivalence
relation, as opposed to the quotient space, used for instance to
define algebraic spaces (for more details and the definition of
equivalence relation see appendix A).

\subsection{Properties of Algebraic Stacks}
\label{subsproperties}

So far we have only defined scheme-theoretic properties for
representable stacks and morphisms. We can define some properties for
arbitrary algebraic stacks (and morphisms among them) using the atlas.

Let ``P'' be a property of schemes, local in nature for the smooth
(resp. \'etale) topology. For example: regular, normal, reduced, of
characteristic $p$,... Then we say that an Artin
(resp. Deligne-Mumford) stack has ``P'' iff the atlas has ``P''
(\cite[p.25]{La}, \cite[p.100]{DM}).

Let ``P'' be a property of morphisms of schemes, local on source and
target for the smooth (resp. \'etale) topology, i.e. for any
commutative diagram
$$
\xymatrix{
{X'} \ar[r]^{p} \ar[dr]_{f''}& {Y'\times_Y X} \ar[r]^{g'} \ar[d]_{f'} 
& {X} \ar[d]^{f} \\
& {Y'} \ar[r]^{g} & {Y} }
$$
with $p$ and $g$ smooth (resp. \'etale) and surjective, $f$ has ``P''
iff $f''$ has ``P''. For example: flat, smooth, locally of finite
type,... For the \'etale topology we also have: \'etale,
unramified,... Then if $f:\SX \to \SY$ is a morphism of Artin
(resp. Deligne-Mumford) stacks, we say that $f$ has ``P'' iff for one
(and then for all) commutative diagram of stacks
$$
\xymatrix{
{X'} \ar[r]^{p} \ar[dr]_{f''}& {Y'\times_Y \SX} \ar[r]^{g'} \ar[d]_{f'} 
& {\SX} \ar[d]^{f} \\
& {Y'} \ar[r]^{g} & {\SY} }
$$
where $X'$, $Y'$ are schemes and $p$, $g$  are smooth (resp. \'etale)
and surjective, $f''$ has ``P'' (\cite[pp. 27-29]{La}).

For Deligne-Mumford stacks it is enough to find a commutative diagram 
$$
\xymatrix{
{X'} \ar[r]^{p} \ar[d]_{f''}& {\SX} \ar[d]^{f} \\
 {Y'} \ar[r]^{g} & {\SY} }
$$
where $p$ and $g$ are \'etale and surjective and $f''$ has ``P''.
Then it follows that $f$ has ``P'' (\cite[p. 100]{DM}).

Other notions are defined as follows.

\begin{definition}[Substack]
\label{substack}
\textup{\cite[def 2.5]{La}, \cite[p.102]{DM}}. 
A stack $\SE$ is a substack of $\SF$ if it is a full subcategory of
$\SF$ and
\begin{enumerate}
\item If an object $X$ of $\SF$ is in $\SE$, then all isomorphic
objects are also in $\SE$.

\item For all morphisms of schemes $f:U\to V$, if $X$ is in $\SE(V)$,
then $f^* X$ is in $\SE(U)$.

\item Let $\{U_i \to U\}$ be a cover of $U$ in the site
$(\Sch/S)$. Then $X$ is in $\SE$ iff $X|_i$ is in $\SE$ for all $i$.
\end{enumerate}
\end{definition}

\begin{definition}
\textup{\cite[def 2.13]{La}}.
A substack $\SE$ of $\SF$ is called open (resp. closed, resp. locally
closed) if the inclusion morphism $\SE \to \SF$ is
\textbf{representable} and it is an open immersion (resp. closed
immersion, resp. locally closed immersion).
\end{definition}

\begin{definition}[Irreducibility]
\textup{\cite[def 3.10]{La}, \cite[p.102]{DM}}.  
An algebraic stack $\SF$ is irreducible if it is not the union of two
distinct and nonempty proper closed substacks.
\end{definition}

\begin{definition}[Separatedness]
\textup{\cite[def 3.17]{La}, \cite[def 4.7]{DM}}.
An algebraic stack $\SF$ is separable, if the (representable) diagonal
morphism $\Delta_\SF$ is universally closed (and hence proper, because
it is automatically separable and of finite type).

A morphism $f:\SF \to \SG$ of algebraic stacks is separable if for all $U
\to \SF$ with $U$ affine, $U\times_\SG \SF$ is a separable (algebraic) stack.
\end{definition}

For Deligne-Mumford stacks, $\Delta_\SF$ is universally closed iff it
is finite.
There is a valuative criterion of separatedness, similar to the
criterion for schemes. Recall that by Yoneda lemma (lemma
\ref{yoneda}), a morphism $f:U\to \SF$ between a scheme and a stack is
equivalent to an object in $\SF(U)$. Then we will say that $\alpha$ is
an isomorphism between two morphisms $f_1,f_2:U\to \SF$ when $\alpha$
is an isomorphism between the corresponding objects of $\SF(U)$.

\begin{proposition}[Valuative criterion of separatedness (stacks)] 
\textup{\cite[prop
3.19]{La}, \cite[thm 4.18]{DM}}.
An algebraic stack $\SF$ is separated  (over $S$) if and
only if the following holds. Let $A$ be a valuation ring with fraction
field $K$.
Let $g^{}_1:\Spec A\to \SF$ and
$g^{}_2:\Spec A \to \SF$ be two morphisms such that:
\begin{enumerate}
\item $f_{p^{}_\SF}\circ g^{}_1= f_{p^{}_\SF}\circ g^{}_2$.
\item There exists an isomorphism $\alpha: g^{}_1|_{\Spec K} 
\to g^{}_2|_{\Spec K}$.
\end{enumerate}
$$
\xymatrix{
 & & {\SF} \ar[d]^{p^{}_{\SF}} \\
{\Spec K} \ar@(u,l)[rru]   \ar[r]^{i} &
{\Spec A} \ar@<0.5ex>[ru]^{g^{}_1} \ar@<-0.5ex>[ru]_{g^{}_2} \ar[r] 
& S
}
$$
then there exists an isomorphism (in fact unique) $\tilde\alpha: 
g^{}_1\to g^{}_2$ that
extends $\alpha$, i.e. $\tilde\alpha|_{\Spec K}=\alpha$. 
\end{proposition}

\begin{remark}
\label{dvr}
\textup{
It is enough to consider complete valuation rings $A$ with
algebraically  closed residue field \cite[3.20.1]{La}. If furthermore
$S$ is locally Noetherian and $\SF$ is locally is finite type, it is
enough to consider discrete valuation rings $A$ \cite[3.20.2]{La}.}
\end{remark}

\begin{example}
\textup{The stack $BG$ won't be separated if $G$ is not proper over
$S$ \cite[3.20.3]{La}, and since we assumed $G$ to be affine, this
won't happen if it is not finite.}

\textup{In general the moduli stack of vector bundles $\BBund$ 
is not separated. It is easy to find
families of vector bundles that contradict the criterion.}

\textup{The stack of stable curves $\overline\SM_g$ is separated 
\cite[prop 5.1]{DM}.}
\end{example}

The criterion for morphisms is more involved  because we are
working with stacks and we have to keep track of the isomorphisms.

\begin{proposition}[Valuative criterion of separatedness (morphisms)] 
\textup{\cite[prop 3.19]{La}}
A morphism of algebraic stacks $f:\SF \to \SG$ is separated if and
only if the following holds. Let $A$ be a valuation ring with fraction
field $K$. Let $g^{}_1:\Spec A\to \SF$ and
$g^{}_2:\Spec A \to \SF$ be two morphisms such that:
\begin{enumerate}
\item There exists an isomorphism $\beta: f\circ g^{}_1\to f\circ g^{}_2$.
\item There exists an isomorphism $\alpha: g^{}_1|_{\Spec K} 
\to g^{}_2|_{\Spec K}$.
\item $f(\alpha)=\beta|_{\Spec K}$.
\end{enumerate}
then there exists an isomorphism (in fact unique) $\tilde\alpha: 
g^{}_1\to g^{}_2$ that
extends $\alpha$, i.e. $\tilde\alpha|_{\Spec K}=\alpha$ and
$f(\tilde\alpha)=\beta$. 
\end{proposition}

Remark \ref{dvr} is also true in this case.

\begin{definition}
\textup{\cite[def 3.21]{La}, \cite[def 4.11]{DM}}.
An algebraic stack $\SF$ is proper (over $S$) if it is separated and of
finite type, and if there is a scheme $X$ proper over $S$ and a
(representable) surjective morphism $X\to \SF$.

A morphism $\SF\to \SG$ is proper if for any affine scheme $U$ and
morphism $U\to \SG$, the fiber product $U\times_\SG \SF$ is proper
over $U$.
\end{definition}

For properness we only have a satisfactory criterion for stacks (see
\cite[prop 3.23 and conj 3.25]{La} for a generalization for morphisms).

\begin{proposition}[Valuative criterion of properness]
\textup{\cite[prop 3.23]{La}, \cite[thm 4.19]{DM}}.
Let $\SF$ be a separated algebraic stack (over $S$). It is proper
(over $S$) if and only if the following condition holds.
Let $A$ be a valuation ring with fraction
field $K$. 
For any commutative diagram
$$
\xymatrix{
 & & {\SF} \ar[d]^{p^{}_\SF} \\
{\Spec K} \ar[r]^{i} \ar[rru]^{g} & {\Spec A} \ar[r] & S }
$$
there exists a finite field extension $K'$ of $K$ such that $g$ extends to
$\Spec(A')$, where $A'$ is the integral closure of $A$ in $K'$. 
$$
\xymatrix{
 & & {\SF} \ar[dd]^{p^{}_\SF} \\
{\Spec K'} \ar[rru]^{g\circ u} \ar[d]_{u} \ar[r] & 
{\Spec A'} \ar[d] \ar@{-->}[ru] \\
{\Spec K} \ar[r]^{i}  & {\Spec A} \ar[r] & S }
$$
\end{proposition}

\begin{example}[Stable curves]
\textup{The Deligne-Mumford stack of stable curves  
$\overline\SM_g$ is proper 
\cite[thm 5.2]{DM}}.
\end{example}

\subsection{Points and dimension}
\label{subspoints}

We will introduce the concept of point of an algebraic
stack and dimension of a stack at a point. The reference for this is
\cite[chapter 5]{La}.

\begin{definition}
Let $\SF$ be an algebraic stack over $S$. The set of points of $\SF$ 
is the set of equivalence classes of pairs $(K,x)$, with $K$ a field
over $S$ (i.e. a field with a morphism of schemes $\Spec K \to S$)
and $x:\Spec K \to \SF$ a morphism of stacks. Two pairs $(K',x')$ and
$(K'',x'')$ are equivalent if there is a field $K$ extension of $K'$
and $K''$ and a commutative diagram
$$
\xymatrix{
{\Spec K} \ar[r] \ar[d] & {\Spec K'} \ar[d]^{x'} \\
{\Spec K''} \ar[r]^{x''} & \SF
}
$$
Given a morphism $\SF \to \SG$ of algebraic stacks and a point of
$\SF$, we define the image of that point in $\SG$ by composition.
\end{definition}

Every point of an algebraic stack is the image of a point of an
atlas. To see this, given a point represented by $\Spec K \to \SF$ 
and an atlas
$X\to \SF$, take any point $\Spec K' \to X\times_\SF \Spec K$. The
image of this point in $X$ maps to the given point.

To define the concept of dimension, recall that if $X$ and $Y$ are
locally Noetherian schemes and $f:X\to Y$ is flat, then for any point
$x\in X$ we have
$$
\dim_x(X)= \dim_x(f) + \dim_{f(x)}(Y),
$$
with $\dim_x(f)=\dim_x(X_{f(x)})$, where $X_y$ is the fiber of $f$
over $y$. 

\begin{definition}
Let $f:\SF\to \SG$ be a representable morphism, locally of finite type,
between two algebraic spaces. Let $\xi$ be a point of $\SF$. Let $Y$
be an atlas of $\SG$
Take a point $x$ in the algebraic space $Y\times_\SG \SF$ that maps to
$\xi$, 
$$
\xymatrix{
{Y\times_\SG \SF} \ar[r] \ar[d]_{\tilde f} & \SF \ar[d]^{f} \\
{Y} \ar[r] & \SG
}
$$
and define the dimension of the morphism $f$ at the point $\xi$ as
$$
\dim_\xi(f)=\dim_x(\tilde f).
$$
\end{definition}

It can be shown that this definition is independent of the choices
made.

\begin{definition}
Let $\SF$ be a locally Noetherian algebraic stack and $\xi$ a point of
$\SF$. Let $u: X\to \SF$ be an atlas, and $x$ a point of $X$ mapping
to $\xi$. We define the dimension of $\SF$ at the point $\xi$ as
$$
\dim_\xi(\SF)=\dim_x(X)-\dim_x(u).
$$
The dimension of $\SF$ is defined as
$$
\dim(\SF)=\operatorname{Sup}_{\xi} (\dim_\xi(\SF)).
$$
\end{definition}

Again, this is independent of the choices made. 

\begin{example}[Quotient by group action]
\textup{
Let $X$ be a smooth scheme of dimension $\dim(X)$ and $G$ a smooth
group of dimension $\dim(G)$ acting on $X$. Let $[X/G]$ be the
quotient stack defined in example \ref{quotient}. Using the atlas 
defined in example 
\ref{atlasquotient}, we see that 
$$
\dim[X/G]=\dim(X)-\dim(G).
$$
Note that we haven't made any assumption on the action. In particular,
the action could be trivial. The dimension of an algebraic stack can
then be negative. For instance, the dimension of the classifying stack
$BG$ defined in example \ref{quotient} has dimension $\dim(BG)=-\dim(G)$.
}
\end{example}

\subsection{Quasi-coherent sheaves on stacks}
\label{subssheaves}

\begin{definition}
\textup{\cite[def 7.18]{Vi}, \cite[def 6.11, prop 6.16]{La}.}
A quasi-coherent sheaf $\calS$ on an algebraic stack $\SF$ is the
following set of data:
\begin{enumerate}

\item For each morphism $X\to \SF$ where $X$ is a scheme, a
quasi-coherent sheaf $\calS_X$ on $X$.

\item For each commutative diagram
$$
\xymatrix{
{X} \ar[r]^f \ar[dr] & {Y} \ar[d] \\
 & {\SF}
}
$$
an isomorphism $\varphi^{}_f: \calS_X \stackrel{\isom}{\too} f^*\calS_Y$,
satisfying the cocycle condition, i.e. for any commutative diagram
\begin{eqnarray}
\label{sheaf2}
\xymatrix{
{X} \ar[r]^{f} \ar[dr] & {Y} \ar[d] \ar[r]^{g}& {Z} \ar[dl] \\
 & {\SF}
}
\end{eqnarray}
we have $\varphi^{}_{g\circ f} = \varphi^{}_f \circ f^* \varphi^{}_g$.
\end{enumerate}

We say that $\calS$ is coherent (resp. finite type, finite
presentation, locally free) if $\calS_X$ is coherent (resp. finite
type, finite presentation, locally free) for all $X$.

A morphism of quasi-coherent sheaves $h:\calS \to \calS'$ is a
collection of morphisms of sheaves $h^{}_X:\calS^{}_X \to \calS'_X$
compatible with the isomorphisms $\varphi$

\end{definition}

\begin{remark}\textup{
Since a sheaf on a scheme can be obtained by glueing the restriction
to an affine cover, it is enough to
consider affine schemes.}
\end{remark}

\begin{example}[Structure sheaf]
\textup{
Let $\SF$ be an algebraic stack. The structure sheaf $\SO_\SF$ is
defined by taking $(\SO_\SF)_X=\SO_X$.}
\end{example}

\begin{example}[Sheaf of differentials]
\textup{
Let $\SF$ be a Deligne-Mumford stack. To define the sheaf of
differentials $\Omega_\SF$, if $U\to \SF$ is an \'etale morphism we
set $(\Omega_\SF)_U=\Omega_U$, the sheaf of differentials of the
scheme $U$. If $V \to \SF$ is another \'etale
morphism and we have a commutative diagram
$$
\xymatrix{
{U} \ar[r]^f \ar[dr] & {V} \ar[d] \\
 & {\SF}
}
$$
then $f$ has to be \'etale,  there is a canonical isomorphism
$\varphi^{}_f :\Omega_{U/S} \to f^* \Omega_{V/S}$, and these canonical
isomorphisms satisfy the cocycle condition.}

\textup{Once we have defined $(\Omega_\SF)_U$ for \'etale morphisms
$U\to \SF$, we can extend the definition for any morphism $X\to \SF$
with $X$ an arbitrary scheme as follows: take an (\'etale) atlas
$U=\coprod U_i \to \SF$. Consider the composition morphism
$$
X\times_\SF U \stackrel{p_2}{\too} U \too \SF,
$$
and define $(\Omega_\SF)_{X\times_\SF U}=p^*_2\Omega_U$. The cocycle
condition for $\Omega_{U_i}$ and \'etale descent implies that 
$(\Omega_\SF)_{X\times_\SF U}$ descends to give a sheaf
$(\Omega_{\SF})_X$ on $X$. It is easy to check that this doesn't depend
on the atlas $U$ used, and that given a commutative diagram like
(\ref{sheaf2}), there are canonical isomorphisms $\varphi$ satisfying
the cocycle condition.}

\end{example}

\begin{example}[Universal vector bundle]
\textup{
Let $\BBund$ be the moduli stack of vector bundles on a scheme $X$
defined in \ref{bbund}. The universal vector bundle $V$ on $\BBund
\times X$ is defined as follows:}

\textup{
Let $B$ be a scheme and $f=(f_1,f_2): B \to \BBund \times X$ a
morphism. By lemma
\ref{yoneda}, the morphism $f_1:B\to \BBund$ is equivalent to a vector
bundle $W$ on $B \times X$. We define $V_B$ as ${\tilde f}^*W$,
where ${\tilde f}=(\id_B, f_2):B\to B\times X$.
Let
$$
\xymatrix{
{B'} \ar[r]^g \ar[dr]_{f'} & {B} \ar[d]^{f} \\
 & {\BBund \times X}
}
$$
be a commutative diagram. Recall that this means that there is an
isomorphism $\alpha:f \circ g \to f'$, and looking at the projection
to $\BBund$ we have an isomorphism $\alpha^{}_1:f^{}_1\circ g \to f'_1$. 
Using lemma \ref{yoneda}, $f^{}_1\circ g$ and $f'_1$ correspond
respectively to the vector
bundles $(g\times \id_X)^*W$ and $W'$ on $B'\times X$, and (again by
lemma \ref{yoneda}) $\alpha^{}_1$
gives an isomorphism between them. It is easy to check that these
isomorphisms satisfy the cocycle condition for diagrams of the form
(\ref{sheaf2}).
}

\end{example}

\section{Vector bundles: moduli stack vs. moduli scheme}
\label{sectionversus}

In this section we will compare, in the context of vector bundles, the
new approach of stacks versus the standard approach of moduli schemes
via geometric invariant theory (GIT).

Fix a scheme $X$, a positive integer $r$ and 
classes $c_i\in
H^{2i}(X)$. 
All vector bundles over $X$ in this section will have rank $r$ and
Chern classes $c_i$. We will also consider vector bundles
on products $B\times X$ where $B$ is a scheme. We will always assume
that these vectors bundles are flat over $B$, and that the restriction
to the slices $\{p\}\times X$ are vector bundles with rank $r$ and
Chern classes $c_i$. Fix also a polarization on $X$. All references to
stability or semistability of vector bundles will mean Gieseker
stability with respect to this fixed polarization.

Recall that the functor $\FBund^{s}$ (resp. $\FBund^{ss}$) 
is the functor from
$(Sch/S)$ to $(Sets)$ that for each scheme $B$ gives the set of
\textit{equivalence} classes of vector bundles over $B\times X$,
flat over $B$ and such that the restrictions $V|_b$ to the slices
$p\times X$ are stable (resp. semistable) vector bundles with fixed
rank and Chern classes, where two vector bundles $V$ and $V'$ on 
$B\times X$ are considered \textit{equivalent} 
if there is a line bundle $L$ on
$B$ such that $V$ is isomorphic to $V'\otimes p^*_B L$.

\begin{theorem}
There are schemes $\Bund^{s}$ and $\Bund^{ss}$, called moduli schemes,
corepresenting the
functors ${\FBund}^{s}$ and ${\FBund}^{ss}$.
\end{theorem}

The moduli scheme $\Bund^{ss}$ is constructed using the Quot schemes
introduced in example \ref{quotconstruction} (for a detailed
exposition of the construction, see \cite{HL}). Since the set of
\textit{semistable} vector bundles is bounded, we can choose once and
for all $N$ and $m$ (depending only on the Chern classes and rank)
with the property that for any semistable vector bundle $V$ there is a
point in $R=R_{N,m}$ whose corresponding quotient is isomorphic to $V$.

The scheme $R$ parametrizes vector bundles $V$ on $X$ 
together with a basis of $H^0(V(m))$ (up to
multiplication by scalar). Recall that $N=h^0(V(m))$. 
There is an action of $\gl$ on $R$,
corresponding to change of basis
but since two basis that only
differ by a scalar give the same point on $R$, this $\gl$ action
factors through $\pgl$. Then the moduli scheme $\Bund^{ss}$
is defined as the GIT quotient $R \gitq \pgl$.

The closed points of $\Bund^{ss}$ correspond to S-equivalence
classes of vector bundles, so if there is a strictly semistable vector
bundle, the functor ${\FBund}^{ss}$ is not representable.

Now we will compare this scheme with the moduli stack $\BBund$ defined on 
example \ref{bbund}. We will also consider the moduli stack
$\BBund^{s}$ defined in the same way, but with the extra requirement
that the vector bundles should be stable. The moduli stack
$\BBund^{s}$ is a substack (definition \ref{substack}) of $\BBund$. 
The following are some of the
differences between the moduli scheme and the
moduli stack:

\begin{enumerate}

\item
The stack $\BBund$ parametrizes all vector bundles, but the scheme 
$\Bund^{ss}$ only
parametrizes semistable vector bundles.

\item
From the point of view of the scheme $\Bund^{ss}$, we identify two
vector bundles if they are S-equivalent. On the other hand, from the
point of view of the moduli stack, two vector bundles are identified 
only if they are isomorphic.

\item
Let $V$ and $V'$ be two families of vector bundles parametrized by a
scheme $B$, i.e. two vector bundles (flat over $B$) on $B\times X$. 
If there is a line bundle $L$ on $B$ such that $V$ is isomorphic to
$V'\otimes p^*_B L$, then from the point of view of the moduli scheme,
$V$ and $V'$ are identified as being the same family. On the other
hand, from the point of view of the moduli stack, $V$ and $V'$ are
identified only if they are isomorphic as vector bundles on $B\times
X$.

\item 
The subscheme $\Bund^{s}$ corresponding to stable vector bundles is
sometimes representable by a scheme, but the moduli stack $\BBund^{s}$
is never representable by a scheme. To see this, note that any vector
bundle has automorphisms different from the identity (multiplication
by scalars) and apply lemma \ref{nonrepresentable}.

\end{enumerate}

Now we will restrict our attention to stable bundles, i.e. to the
scheme $\Bund^s$ and the stack $\BBund^s$. For stable bundles the
notions of $S$-equivalence and isomorphism coincide, so the points of
$\Bund^s$ correspond to isomorphism classes of vector
bundles. Consider $R^{s}\subset R$, the subscheme corresponding to stable
bundles. There
is a map $\pi :R^s \to \Bund^s=R^s/\pgl$, and $\pi$ is in fact a
principal $\pgl$-bundle (this is a consequence of Luna's \'etale slice
theorem).

\begin{remark}[Universal bundle on moduli scheme]\textup{
The scheme $\Bund^s$ represents the functor $\FBund^s$ if there is a
universal family. Recall that a universal family for this functor is a vector
bundle $E$ on $\Bund^s \times X$ such that the isomorphism class of
$E|_{p\times X}$ is the isomorphism class corresponding to the point
$p\in \Bund^s$, and for any family of vector bundles $V$ on $B\times
X$ there is a morphism $f:B\to \Bund^s$ and a line bundle $L$ on $B$ 
such that $V \otimes p^*_B L$ is isomorphic to $(f\times \id)^*E$.
Note that if $E$ is a universal family, then $E\otimes p^*_{\Bund^s}L$
will also be a universal family for any line bundle $L$ on $\Bund^s$.}

\textup{
The universal bundle for the Quot scheme gives a universal family 
$\wt V$ on $R^s\times X$, but this family doesn't always descend to
give a universal family on the quotient $\Bund^s$.}

\textup{
Let $X\stackrel{G}\too Y$ be a principal $G$-bundle. A vector bundle
$V$ on $X$ descends to $Y$ if the action of $G$ on $X$ can be lifted
to $X$. In  our case, if certain numerical criterion involving $r$ and
$c_i$ is satisfied (if $X$ is a smooth curve this criterion is
$\operatorname{gcd}(r,c_1)=1$), then we can find a line bundle $L$ on
$R^s$ such that the $\pgl$ action on $R^s$ can be lifted to
$\wt V \otimes p^*_{R^s}L$, and then this vector bundle descends to
give a universal family on $\Bund^s \times X$. But in general the best
that we can get is a universal family on an \'etale cover of
$\Bund^s$.}
\end{remark}

Recall from example \ref{atlasquotient} that there is a morphism
$[R^{ss}/\pgl] \to \Bund^{ss}$, and that the morphism
$[R^{s}/\pgl] \to \Bund^{s}$ is an isomorphism of stacks.

\begin{proposition}
\label{versus}

There is a commutative diagram of stacks
$$
\xymatrix{
{[R^{s}/\gl]} \ar[rr]^{q} \ar[d]_{g}^{\simeq}& 
&{[R^{s}/\pgl]} \ar[d]^{h}_{\simeq} \\
{\BBund^{s}} \ar[rr]_{\varphi} &  &{\;\Bund^{s},}
}
$$
where $g$ and $h$ are isomorphisms of stacks, but $q$ and $\varphi$
are not.
If we change ``stable'' by
``semistable'' we still have a commutative diagram, 
but the corresponding morphism $h^{ss}$ is not an
isomorphism of stacks.

\end{proposition}

\begin{proof}
The morphism $\varphi$ is the composition of the natural morphism
$\BBund^{s} \to \FBund^{s}$ (sending each category to the set of
isomorphism classes of objects) and the morphism $\FBund^{s} \to 
\Bund^{s}$ given
by the fact that the scheme $\Bund^{s}=R^s\gitq \pgl$ corepresents the 
functor.

The morphism $h$ was constructed in example \ref{quotient}.

The key ingredient needed to define $g$ is the fact that the $\gl$
action on the Quot scheme lifts to the universal bundle, i.e.
the universal bundle on the Quot scheme has a $\gl$-linearization.
Let
$$
\xymatrix{
{\wt{B}} \ar[r]^{p} \ar[d] & R^{ss} \\
{B}  }
$$
be an object of $[R^{ss}/\gl]$. Since $R^{ss}$ is a subscheme of a
Quot scheme, and this universal bundle has a $\gl$-linearization. 
Let $\wt E$ be the vector
bundle on $\wt B\times X$ defined by the pullback of this universal
bundle. Since $f$ is $\gl$-equivariant, $\wt E$ is also
$\gl$-linearized. Since $\wt B \times X \to B\times X$ is a principal
bundle, the vector bundle $\wt E$ descends to give a vector bundle $E$
on $B\times X$, i.e. an object of $\BBund^{ss}$. Let
$$
\xymatrix{ & & R^{ss}\\
{\wt{B}} \ar[r]_{\phi} \ar[d] \ar[rru]^{f}
& {\wt{B}'} \ar[d] \ar[ru]_{f'} \\
{B} \ar@{=}[r] & {B} }
$$
be a morphism in $[R^{ss}/\gl]$. Consider the vector bundles $\wt E$
and $\wt E'$ defined as before. Since $f'\circ \phi=f$, we get an
isomorphism of $\wt E$ with $(\phi \times \id)^* \wt E'$. Furthermore
this isomorphism is $\gl$-equivariant, and then it descends to give an
isomorphism of the vector bundles $E$ and  $E'$ on $B\times X$, and
we get a morphism in $\BBund^{ss}$.

To prove that this gives an equivalence of categories, we construct a
functor $\overline g$ from $\BBund^{ss}$ to $[R^{ss}/\gl]$.
Given a vector bundle on $B\times X$, let $q:\wt B \to B$ be the
$\gl$-principal bundle associated with the vector bundle $p^*_B E$ on
$B$. Let $\wt E=(q\times \id)^*E$ be the pullback of $E$ to $\wt
B\times X$. It has a canonical $\gl$-linearization because it is
defined as a pullback by a principal $\gl$-bundle. The vector bundle
${p^{}_{\wt B}}_* \wt E$ is canonically isomorphic to the trivial
bundle $\SO^N_{\wt B}$, and this isomorphism is $\gl$-equivariant, so
we get an \textit{equivariant} morphism $\wt B\to R^{ss}$, and hence an object
of $[R^{ss}/\gl]$.

If we have an isomorphism between two vector bundles $E$ and $E'$ on
$B\times X$, it is easy to check that it induces an isomorphism
between the associated objects of $[R^{ss}/\gl]$.

It is easy to check that there are natural isomorphisms of functors
 $g\circ \wt g \isom \id$ and $\wt g\circ g \isom \id$, 
and then $g$ is an equivalence of
categories.

The morphism $q$ is defined using the following lemma, with $G=\gl$,
$H$ the subgroup consisting of scalar multiples of the identity,
$\overline G=\pgl$ and $Y$=$R^{ss}$.

\end{proof}

\begin{lemma}
Let $Y$ be an $S$-scheme and $G$ an affine flat group $S$-scheme, 
acting on $Y$ on the right. Let $H$ be a normal closed 
subgroup of
$G$. Assume that $\overline G=G/H$ is affine. If $H$ acts trivially 
on $Y$, then there is a
morphism of stacks
$$
[Y/G]\too [Y/\overline G].
$$
If $H$ is nontrivial, then this morphism is not faithful, 
so it is not an isomorphism.
\end{lemma}

\begin{proof}
Let 
$$
\xymatrix{
{E} \ar[r]^{f} \ar[d]^{\pi} & Y \\
{B}  }
$$
be an object of $[Y/G]$. There is a scheme $Y/H$ such that $\pi$
factors
$$
E \stackrel{q}\too E/H \stackrel{\pi'}\too B.
$$
To construct $Y/H$, note that there is a local \'etale cover $U_i$ of
$B$ and isomorphisms $\phi_i:\pi^{-1}(U_i)\to U_i\times G$, with transition
functions $\psi_{ij}=\phi^{}_i \circ \phi^{-1}_j$. Since these
isomorphisms are $G$-equivariant, they descend to give isomorphisms
$\overline{\psi}_{ij}:U_j\times G/H \to U_i\times G/H$, and using this
transition functions we get $Y/H$. This construction shows that $\pi'$
is a principal $\overline G$-bundle. 
Furthermore, $q$ is also a
principal $H$-bundle (\cite[example 4.2.4]{HL}), and in particular it
is a categorical quotient.

Since $f$ is $H$-invariant, there is a morphism $\overline f: E/H \to
R$, and this gives an object of $[Y/\overline G]$.

If we have a morphism in $[Y/G]$, given by a morphism $g:E\to E'$ of
principal $G$-bundles over $B$, it is easy to see that it descends
(since $g$ is equivariant) to
a morphism $\overline{g}:E/H \to E'/H$, giving a morphism in
$[Y/\overline G]$.

This morphism is not faithful, since the automorphism
$E\stackrel{\cdot z}{\too} E$ given by
multiplication on the right by a nontrivial element $z\in H$
is sent
to the identity automorphism $E/H \to E/H$, and then $\Hom(E,E)\to
\Hom(E/H,E/H)$ is not injective.

\end{proof}

If $X$ is a smooth curve, then it can be shown that $\BBund$ is a
smooth stack of dimension $r^2(g-1)$, where $r$ is the rank and $g$ is
the genus of $X$. In particular, the open substack $\BBund^{ss}$ is
also smooth of dimension $r^2(g-1)$, but the moduli scheme $\Bund^{ss}$ is
of dimension $r^2(g-1)+1$ and might not be smooth. Proposition
\ref{versus} explains the difference in the dimensions (at least on
the smooth part): we obtain the moduli stack by taking the quotient by
the group $\gl$, of dimension $N^2$, but the moduli scheme is obtained
by a quotient by the group $\pgl$, of dimension $N^2-1$. The moduli
scheme $\Bund^{ss}$ is not smooth in general because in the strictly
semistable part of $R^{ss}$ the action of $\pgl$ is not free. On the
other hand, the smoothness of a stack quotient doesn't depend on the
freeness of the action of the group.

\section{Appendix A: Grothendieck topologies, sheaves and algebraic spaces}
\label{grothendiecktopologies}

The standard reference for Grothendieck topologies is SGA 
(\textit{S\'e\-mi\-naire de G\'eo\-m\'e\-trie Alg\'e\-bri\-que}). For an
introduction see \cite{T} or \cite{MM}. For algebraic spaces, see
\cite{K} or \cite{Ar1}.

An open cover in a
topological space $U$ can be seen as family of morphisms in the category of
topological spaces $f_i:U_i \to U$, with the property that $f_i$ is an open
inclusion and the union of their images is $U$, i.e we are choosing a 
class of morphisms (open inclusions)
in the category of topological spaces.
A Grothendieck topology on an arbitrary category is basically a choice of a
class of morphisms, that play the role of ``open
sets''. A morphism $f:V\to U$ in this class is to be thought of as an
``open set'' in the object $U$. The concept of intersection of open
sets, for instance, can be replaced by the fiber product: the
``intersection'' of $f_1:U_1\to U$ and $f_2:U_2\to U$ is
$f_{12}:U_1\times _U U_2 \to U$.

A category with a Grothendieck topology is called a site. We will
consider two topologies on $(\Sch/S)$.
\medskip

\textbf{fppf topology}. Let $U$ be a scheme. Then a cover of $U$ is a
finite collection of morphisms $\{f_i:U_i\to U\}_{i\in I}$ such that
each $f_i$ is a finitely presented flat morphism (for Noetherian
schemes, this is equivalent to flat and finite type), and $U$ is the
(set theoretic) union of the images of $f_i$. In other words, 
$\coprod U_i \to U$ is \textit{``fid\`element plat de pr\'esentation
finie''}.
\medskip

\textbf{\'Etale topology}. Same definition, but substituting flat by
\'etale.
\medskip

A presheaf of sets on $(\Sch/S)$ is a contravariant functor $F$ from
$(\Sch/S)$ to $(\Sets)$. Choose a topology on $(\Sch/S)$. 
We say that
$F$ is a sheaf (or an $S$-space) with respect to that topology if 
for every cover
$\{f_i:U_i\to U\}_{i\in I}$ in the topology the following two axioms
are satisfied:
\begin{enumerate}
\item \textit{(Mono)} Let $X$ and $Y$ be two elements of $F(U)$. If
$X|_i=Y|_i$ for all $i$, then $X=Y$.

\item \textit{(Glueing)} Let $X_i$ be an object of $F(U_i)$ for each
$i$  such that
$X_i|_{ij}=X_j|_{ij}$, then there exists $X \in F(U)$ such that
$X|_i=X_i$ for each $i$.
\end{enumerate}

We have used the following notation: if $X\in F(U)$, then $X|_i$ is
the element of $F(U_i)$ given by $F(f_i)(X)$, and if $X_i\in F(U_i)$,
then $X_i|_{ij}$ is the element of $F(U_{ij})$ given by
$F(f_{ij,i})(X_i)$ where $f_{ij,i}:U_i\times_U U_j \to U_i$ is the
pullback of $f_j$.

\smallskip
We can define morphisms of $S$-spaces as morphisms of sheaves (natural
transformation of functors with the obvious conditions).
Note that a scheme can be viewed as an $S$-space via its functor of
points, and a morphism between two such $S$-spaces is equivalent to a
scheme morphism between the schemes (by the Yoneda embedding lemma),
then the category of $S$-schemes is a full subcategory of the category of
$S$-spaces.

\medskip
\textbf{Equivalence relation and quotient space}.
An equivalence relation in the category of $S$-spaces consists of two
$S$-spaces $R$ and $U$ and a monomorphism of $S$-spaces
$$
\delta:R \to U \times_S U
$$
such that for all $S$-scheme $B$, the map $\delta(B):R(B)\to
U(B)\times U(B)$ is the graph of an equivalence relation between
sets. A quotient $S$-space for such an equivalence relation is by
definition the sheaf cokernel of the diagram
$$
\xymatrix{
{R} \ar@<0.5ex>[r]^{p_2\circ \delta} 
\ar@<-0.5ex>[r]_{p_1\circ \delta} & {U}}
$$

\begin{definition}
\textup{\cite[0]{La}}.
An $S$-space $F$ is called an algebraic space if it is the quotient 
$S$-space
for an equivalence
relation such that $R$ and $U$ are $S$-schemes, $p_1\circ \delta$,
$p_2\circ \delta$ are \'etale (morphisms of $S$-schemes), and $\delta$
is a quasi-compact morphism (of $S$-schemes).
\end{definition}

Roughly speaking, an algebraic space is a quotient of a scheme by an
\'etale equivalence relation. The following is an equivalent
definition.

\begin{definition}
\textup{\cite[def 1.1]{K}}.
An $S$-space $F$ is called an algebraic space if there 
exists a scheme
$U$ (atlas) and a morphism of $S$-spaces $u:U\to F$ such that
\begin{enumerate}
\item (The morphism $u$ is \'etale) For any $S$-scheme $V$ and
morphism $V \to F$, the (sheaf) fiber
product $U\times_F V$ is representable by a scheme, and the map
$U\times_F V\to V$ is an \'etale morphism of schemes.
\item (Quasi-separatedness) The morphism $U\times_F U \to
U\times_S U$ is quasi-compact.
\end{enumerate}
\end{definition}

We recover the first definition by taking $R=U\times_F U$. Then 
roughly speaking, we can also think of an algebraic space as
``something'' that looks locally in the \'etale topology like an 
affine scheme,
in the same sense that a scheme is something that looks locally in
the Zariski topology like an affine scheme.

Algebraic spaces are used, for instance, to give algebraic structure 
to certain
complex manifolds (for instance Moishezon manifolds) that are not 
schemes, but
can be realized as algebraic spaces.
All smooth algebraic spaces of dimension 1 and
2 are actually schemes. An example of a smooth algebraic space of
dimension 3 that is not a scheme can be found in \cite{H}.

But \'etale topology is useful even if we are only interested in
schemes.
The idea is that the \'etale topology is finer than the Zariski
topology, and in many situations it is ``fine enough'' to do the
analogue of the manipulations that can be done with the analytic
topology of complex manifolds. As an example, consider the affine
complex line $\Spec(\CC[x])$, and take a (closed) point $x_0$
different from $0$. Assume that we want to define the function $\sqrtx$
in a neighborhood of $x_0$. In the analytic topology we only need to
take a neighborhood small enough so that it doesn't contain a loop
that goes around the
origin, then we choose one of the branches (a sign) of the square
root. In the Zariski topology this cannot be done, because all open
sets are too large (have loops going around the origin, so the sign of
the square root will change, and $\sqrtx$ will be multivaluated).
But take the 2:1 \'etale map $V= \Spec(\CC[y,x,x^{-1}]/(y-x^2))  \to 
\Spec(\CC[x])$. 
The function $\sqrtx$ can certainly be defined on $V$, it is just
equal to the function $y$, so it is in this sense that we say that 
the \'etale topology is finer: $V$ is a ``small enough open subset''
because the square root can be defined on it.

\section{Appendix B: 2-categories}

In this section we recall the notions of 2-category and 2-functor.
A 2-category $\mathfrak{C}$ consists of the following data \cite{Hak}:

\begin{enumerate}
\item [(i)] A class of objects $\obc$
\item [(ii)] For each pair $X$, $Y \in \obc$, a category $\Hom(X,Y)$
\item [(iii)] \textit{horizontal composition of 1-morphisms and
2-morphisms}. For each triple $X$, $Y$, $Z \in \obc$, a functor
$$
\mu_{X,Y,Z}:\Hom(X,Y) \times \Hom(Y,Z) \to \Hom (X,Z)
$$
\end{enumerate}
with the following conditions
\begin{enumerate}
\item [(i')] \textit{(Identity 1-morphism)} For each object $X\in \obc$,
there exists an object
$\id_X\in \Hom(X,X)$ such that
$$
\mu_{X,X,Y}(\id_X,\;)=\mu_{X,Y,Y}(\;,\id_Y)=\id_{\Hom(X,Y)},
$$
where $\id_{\Hom(X,Y)}$ is the identity functor on the category
$\Hom(X,Y)$
\item[(ii')] \textit{(Associativity of horizontal compositions)}
For each quadruple $X$, $Y$, $Z$, $T\in \obc$,
$$
\mu_{X,Z,T}\circ (\mu_{X,Y,Z}\times \id_{\Hom(Z,T)})=
\mu_{X,Y,T}\circ (\id_{\Hom(X,Y)}\times\mu_{Y,Z,T})
$$
\end{enumerate}

The example to keep in mind is the 2-category $\mathfrak{Cat}$ of
categories. The objects of $\mathfrak{Cat}$ are categories, and for
each pair $X$, $Y$ of categories, $\Hom(X,Y)$ is the category of
functors between $X$ and $Y$.

Note that the main difference between a 1-category (a usual category)
and a 2-category is that $\Hom(X,Y)$, instead of being a set, is a
category.

Given a 2-category, an object $f$ of the category $\Hom(X,Y)$ is
called a 1-morphisms of
$\fc$, and is represented with a diagram
$$
\xymatrix
{
{\bullet} \ar[r]^f \save[]+<0ex,2.5ex>*{X}\restore   & 
{\bullet}\save[]+<0ex,2.5ex>*{Y}\restore}
$$
and a morphism $\alpha$ of the category $\Hom(X,Y)$ is called a 
2-morphisms of $\fc$, and is represented as
$$
\xymatrix
{
{\bullet} \ar @(ur,ul)[rr]^f_{}="f" \ar @(dr,dl)[rr]_{f'}^{}="fp" 
\save[]+<0ex,2.5ex>*{X}\restore   & &{\bullet}
\save[]+<0ex,2.5ex>*{Y}\restore
\ar @2^{\alpha} "f";"fp"}
$$
Now we will rewrite the axioms of a 2-category using diagrams.
\begin{enumerate}
\item \textit{(Composition of 1-morphisms)} Given a diagram
$$
\xymatrix
{{\bullet} \ar[r]^f \save[]+<0ex,2.5ex>*{X}\restore   & 
{\bullet} \ar[r]^g \save[]+<0ex,2.5ex>*{Y}\restore & 
{\bullet} \save[]+<0ex,2.5ex>*{Z}\restore}
\quad\text{there exist}\quad
\xymatrix
{{\bullet} \ar[r]^{g\circ f} \save[]+<0ex,2.5ex>*{X}\restore   & 
{\bullet}\save[]+<0ex,2.5ex>*{Z}\restore}
$$
(this is (iii) applied to objects) and this composition is
associative: $(h\circ g) \circ f= h\circ (g\circ f)$ (this is (ii')
applied to objects). 

\item \textit{(Identity for 1-morphisms)} For each object $X$ there is a
1-morphism $\id_X$ such that $f\circ \id_Y =\id_X \circ f=f$ (this is
(i')).

\item \label{three} \textit{(Vertical composition of 2-morphisms)} 
Given a diagram
$$
\xymatrix
{{\bullet} \ar @(ur,ul)[rr]^f_{}="f" \ar [rr]|g^{}="g"_{}="g2" 
\ar @(dr,dl)[rr]_h^{}="h"
\save[]+<0ex,2.5ex>*{X}\restore   & &{\bullet}
\save[]+<0ex,2.5ex>*{Y}\restore
\ar @2^{\alpha} "f";"g"
\ar @2^{\beta} "g2";"h"}
\quad\text{there exists}\quad
\xymatrix
{
{\bullet} \ar @(ur,ul)[rr]^f_{}="f" \ar @(dr,dl)[rr]_h^{}="g" 
\save[]+<0ex,2.5ex>*{X}\restore   & &{\bullet}
\save[]+<0ex,2.5ex>*{Y}\restore
\ar @2^{\beta\circ\alpha} "f";"g"}
$$
and this composition is associative $(\gamma\circ\beta)\circ\alpha =
\gamma\circ(\beta\circ\alpha)$.

\item \textit{(Horizontal composition of 2-morphisms)} Given a diagram
$$
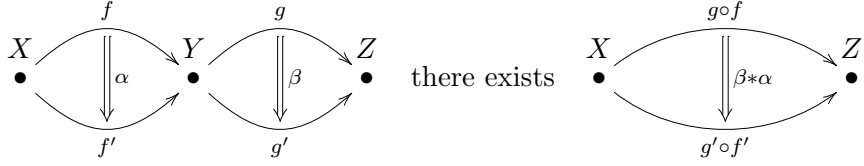
\xymatrix
{
{\bullet} \ar @(ur,ul)[rr]^f_{}="f" \ar @(dr,dl)[rr]_{f'}^{}="fp" 
\save[]+<0ex,2.5ex>*{X}\restore   & &{\bullet}
\save[]+<0ex,2.5ex>*{Y}\restore
 \ar @(ur,ul)[rr]^{g}_{}="g" \ar @(dr,dl)[rr]_{g'}^{}="gp" 
  & &{\bullet}
\save[]+<0ex,2.5ex>*{Z}\restore
\ar @2^{\alpha} "f";"fp"
\ar @2^{\beta} "g";"gp"}
\quad\text{there exists}\quad
\xymatrix
{
{\bullet} \ar @(ur,ul)[rrr]^{g\circ f}_{}="gf" \ar @(dr,dl)[rrr]
_{g'\circ f'}^{}="gpfp" 
\save[]+<0ex,2.5ex>*{X}\restore  & & &{\bullet}
\save[]+<0ex,2.5ex>*{Z}\restore
\ar @2^{\beta\ast\alpha} "gf";"gpfp"}
$$
(this is (iii) applied to morphisms) and it is associative $(\gamma\ast
\beta)\ast\alpha=\gamma\ast(\beta\ast\alpha)$ (this is (ii') applied
to morphisms).

\item \textit{(Identity for 2-morphisms)} For every 1-morphism $f$ there
is a 2-morphism $\id_f$ such that $\alpha\circ\id_g=\id_f\circ\alpha=
\alpha$ (this and item \ref{three} are (ii)). We have $\id_g
\ast \id_f=\id_{g\circ f}$ (this means that $\mu_{X,Y,Z}$ respects the
identity).

\item \textit{(Compatibility between horizontal and vertical
composition of 2-morphisms)} Given a diagram
$$
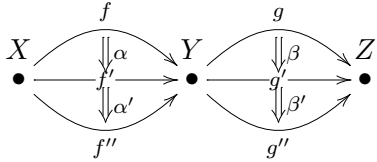
\xymatrix
{{\bullet} \ar @(ur,ul)[rr]^f_{}="f" \ar [rr]|{f'}^{}="f1"_{}="f2" 
\ar @(dr,dl)[rr]_{f''}^{}="fpp"
\save[]+<0ex,2.5ex>*{X}\restore   & &
{\bullet} \ar @(ur,ul)[rr]^g_{}="g" \ar [rr]|{g'}^{}="g1"_{}="g2" 
\ar @(dr,dl)[rr]_{g''}^{}="gpp"
\save[]+<0ex,2.5ex>*{Y}\restore   & &{\bullet}
\save[]+<0ex,2.5ex>*{Z}\restore
\ar @2^{\alpha} "f";"f1"
\ar @2^{\alpha'} "f2";"fpp"
\ar @2^{\beta} "g";"g1"
\ar @2^{\beta'} "g2";"gpp"}
$$
then $(\beta'\circ \beta)\ast(\alpha'\circ \alpha)=(\beta'\ast\alpha')
\circ(\beta\ast\alpha)$ (this is (iii) applied to morphisms).
\end{enumerate}
Two objects $X$ and $Y$ of a 2-category are called equivalent if 
there exist two 1-morphisms
$f:X\to Y$, $g:Y\to X$ and two 2-isomorphisms (invertible
2-morphism) $\alpha:g\circ f \to
\id_X$ and $\beta:f\circ g \to \id_Y$.

A commutative diagram of 1-morphisms in a 2-category is a diagram
$$
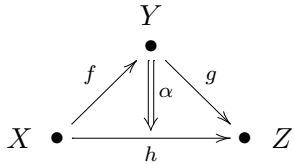
\xymatrix{
 &  {\bullet} \ar[rd]^g \save[]+<0ex,2.5ex>*{Y}\restore
    \ar @2[d]^{\alpha}  \\
{\bullet} \ar[ru]^f \ar[rr]_{h}
   \save[]-<3ex,0ex>*{X}\restore & &
{\bullet} \save[]+<3ex,0ex>*{Z}\restore}
$$
such that $\alpha:g\circ f \to h$ is a 2-isomorphisms.

\begin{remark}
\textup{
Since 2-functors only respect composition of 1-functors up to a
2-isomorphism (condition 3), sometimes they are called pseudofunctors
or lax functors.}
\end{remark}

\begin{remark}
\textup{
Note that we don't require $g\circ f=h$ to say that the diagram is
commutative, but just require that there is a 2-isomorphisms between
them. This is the reason why 2-categories are used to describe
stacks.}
\end{remark}

On the other hand, a diagram of 2-morphisms will be called
commutative only if the compositions are actually equal.
Now we will define the concept of covariant 2-functor (a contravariant
2-functor is defined in a similar way).

A covariant 2-functor $F$ between two 2-categories $\fc$ and $\fcp$ is a law
that for each object $X$ in $\fc$ gives an
object $F(X)$ in $\fcp$. For each 1-morphism $f:X\to Y$ in $\fc$ gives
a 1-morphism $F(f):F(X)\to F(Y)$ in $\fcp$, and for each 2-morphism 
$\alpha:f\Rightarrow g$ in $\fc$ gives a 2-morphism
$F(\alpha):F(f)\Rightarrow F(g)$ in $\fcp$, such that
\begin{enumerate}
\item \textit{(Respects identity 1-morphism)} $F(\id_X)=\id_{F(X)}$.

\item \textit{(Respects identity 2-morphism)} $F(\id_f)=\id_{F(f)}$.

\item \label{twoisom} \textit{(Respects composition of 1-morphism up to a
2-isomorphism)} For every diagram
$$
\xymatrix
{{\bullet} \ar[r]^f \save[]+<0ex,2.5ex>*{X}\restore   & 
{\bullet} \ar[r]^g \save[]+<0ex,2.5ex>*{Y}\restore & 
{\bullet} \save[]+<0ex,2.5ex>*{Z}\restore}
$$
there exists a 2-isomorphism $\epsilon_{g,f}:F(g)\circ F(f) \to
F(g\circ f)$
$$
\xymatrix{
 &  {\bullet} \ar[rd]^{F(g)} \save[]+<0ex,2.5ex>*{F(Y)}\restore
    \ar @2[d]^{\epsilon_{g,f}}  \\
{\bullet} \ar[ru]^{F(f)} \ar[rr]_{F(g\circ f)}
   \save[]-<3ex,0ex>*{F(X)}\restore & &
{\bullet} \save[]+<3ex,0ex>*{F(Z)}\restore}
$$
\begin{enumerate}
\item $\epsilon_{f,\id_X}=\epsilon_{\id_Y,f}=\id_{F(f)}$

\item $\epsilon$ \textit{ is associative}. The following
diagram is commutative
$$
\xymatrix
{F(h)\circ F(g)\circ F(f) \ar@2[rr]^{\epsilon_{h,g} \times \id}
                          \ar@2[d]_{\id \times \epsilon_{g,f}} & &
F(h\circ g)\circ F(f) \ar@2[d]^{\epsilon_{h\circ g,f}} \\
F(h)\circ F(g\circ f) \ar@2[rr]^{\epsilon_{h,g\circ f}} & &
F(h\circ g\circ f)}
$$
\end{enumerate}

\item \textit{(Respects vertical composition of 2-morphisms)}
For every pair of 2-morphisms $\alpha:f \to f'$, $\beta:g \to g'$, 
we have $F(\beta\circ
\alpha)=F(\beta)\circ F(\alpha)$.

\item \label{last} \textit{(Respects horizontal composition of 2-morphisms)} 
For every pair of 2-morphisms $\alpha:f \to f'$, $\beta:g \to g'$, 
the following diagram commutes
$$
\xymatrix
{F(g)\circ F(f) \ar@2[rr]^{F(\beta)\ast F(\alpha)}
                          \ar@2[d]_{\epsilon_{g,f}} & &
F(g')\circ F(f') \ar@2[d]^{\epsilon_{g',f'}} \\
F(g\circ f) \ar@2[rr]^{F(\beta\ast\alpha)} & &
F(g'\circ f')}
$$
\end{enumerate}
By a slight abuse of language, condition \ref{last} is usually written
as $F(\beta)\ast F(\alpha)=F(\beta\ast \alpha)$. Note that strictly
speaking this equality doesn't make sense, because the sources (and
the targets) don't coincide, but if we chose once and for all the
2-isomorphisms $\epsilon$ of condition \ref{twoisom}, then there is a 
unique way
of making sense of this equality.

\begin{remark}
\label{B2}
\textup{
In the applications to stacks, the isomorphism $\epsilon_{g,f}$ of
item \ref{twoisom} is canonically defined, and by abuse of language we will
say that $F(g)\circ F(f)= F(g\circ f)$, instead of saying that they
are isomorphic.}
\end{remark}

Given a 1-category $C$ (a usual category), we can define a 2-category:
we just have to make the set $\Hom(X,Y)$ into a category, and we do
this just by defining the unit morphisms for each element.

On the other hand, given a 2-category $\fc$ there are two ways of
defining a 1-category. We have to make each category $\Hom(X,Y)$ into a
set. The naive way is just to take the set of objects of $\Hom(X,Y)$,
and then we obtain what is called the underlying category of $\fc$ 
(see \cite{Hak}). This has the problem that a 2-functor $F:\fc \to \fcp$
is not in general a functor of the underlying categories (because in
item \ref{twoisom} we only require the composition of 1-morphisms to 
be respected up to 2-isomorphism).

The best way of constructing a 1-category from a 2-category is to
define the set of morphisms between the objects $X$ and $Y$ as the 
set of isomorphism classes of objects of $\Hom(X,Y)$: two objects
$f$ and $g$ of $\Hom(X,Y)$ are isomorphic if there exists a
2-isomorphism $\alpha:f \Rightarrow g$ between them. We call the category
obtained in this way the 1-category associated to $\fc$. Note that a
2-functor between 2-categories then becomes a functor between the
associated 1-categories.

\bigskip

\bigskip

\textbf{Acknowledgments.}
This article is based on a series of lectures that I gave in February
1999 in the Geometric Langlands programme seminar of the Tata Institute 
of Fundamental Research.
First of all, I would like to thank N. Nitsure for proposing
me to give these lectures.
Most of my understanding on stacks comes from
conversations with N. Nitsure and C. Sorger.

I would also like to thank T.R. Ramadas for encouraging me 
to write these notes, and the participants in the
seminar in TIFR for their active participation, interest, 
questions and comments.
In ICTP (Trieste) I gave two informal talks in August 1999
on this subject, and the comments of the participants, specially L.
Brambila-Paz and Y.I. Holla, helped to 
remove mistakes and improve the original notes.

This work was supported by a postdoctoral fellowship of
Ministerio de Educaci\'on y Cultura (Spain).

\end{document}